\documentclass[]{interact}

\usepackage[caption=false]{subfig}
\usepackage[numbers,sort&compress]{natbib}

\usepackage{color}
\usepackage{subfiles}
\usepackage{graphicx}
\usepackage{hyperref}
\usepackage{stmaryrd}
\usepackage{amsmath,amssymb,amsfonts,amsthm,textcomp}
\usepackage{mathtools}
\usepackage{setspace}
\usepackage{nicefrac}
\usepackage{algpseudocode}
\usepackage{algorithm}
\usepackage{booktabs}
\usepackage{bbm}

\bibpunct[, ]{[}{]}{,}{n}{,}{,}

\newfont{\tenbfsl}{cmbxti9 scaled 1200}
\newfont{\tenbbb}{msbm10}
\newfont{\svnbbb}{msbm8}

\newcommand{\Expc}[2]{\mathbb{E}\left[ \left. #1 \right| #2 \right]}

\newcommand{\bs}[1]{\boldsymbol{#1}}
\newcommand{\cl}[1]{\mathcal{#1}}

\newcommand{\bb}[1]{\mathbb{#1}}

\newcommand{\grad}[1]{\nabla F(#1)}
\newcommand{\hess}[1]{\nabla^2 F(#1)}

\newcommand{\xik}{\bs{\xi}_k}
\newcommand{\xikp}{\bs{\xi}_{k+1}}

\newcommand{\bigO}[1]{\cl{O}\left( #1 \right)}
\newcommand{\Hk}{\hat{\bs{H}_k}}
\newcommand{\Bk}{\hat{\bs{B}_k}}
\newcommand{\dimxi}{d_{\bs{\xi}}}

\newcommand{\norm}[1]{\left\lVert {#1} \right\rVert}

\newcommand{\trans}{\scriptscriptstyle\mskip-1mu\top\mskip-2mu}

\newcommand{\sym}{\hbox{\textrm{sym}}\mskip2mu}

\newtheorem{thm}{Theorem}

\theoremstyle{remark}

\newtheorem{ass}{Assumption}

\newcommand{\rset}{\bb{R}}

\begin{document}

\title{Approximating Hessian matrices using Bayesian inference:
a new approach for quasi-Newton methods in stochastic optimization}
\author{
    \name{
        Andr\'{e} Carlon \textsuperscript{a}
        \thanks{CONTACT A.~N. Author. Email: agcarlon@gmail.com}
        Luis Espath \textsuperscript{b}
        Ra\'{u}l Tempone \textsuperscript{a, c, d}
    }
    \affil{
        \textsuperscript{a}
        King Abdullah University of Science \& Technology (KAUST), Computer, Electrical and Mathematical Sciences \& Engineering Division (CEMSE), Thuwal 23955-6900, Saudi Arabia;
        \textsuperscript{b}
        School of Mathematical Sciences, University of Nottingham, Nottingham, NG7 2RD, United Kingdom;
        \textsuperscript{c}
        Department of Mathematics, RWTH Aachen University, Geb\"{a}ude-1953 1.OG, Pontdriesch 14-16, 161, 52062 Aachen, Germany;
        \textsuperscript{d}
        Alexander von Humboldt Professor in Mathematics for Uncertainty Quantification, RWTH Aachen University, Germany
    }
    }


\begin{abstract}
Using quasi-Newton methods in stochastic optimization is not a trivial task given the difficulty of extracting curvature information from the noisy gradients.
Moreover, pre-conditioning noisy gradient observations tend to amplify the noise.
We propose a Bayesian approach to obtain a Hessian matrix approximation for stochastic optimization that minimizes the secant equations residue while retaining the extreme eigenvalues between a specified range.
Thus, the proposed approach assists stochastic gradient descent to converge to local minima without augmenting gradient noise.
We propose maximizing the log posterior using the Newton-CG method.
Numerical results on a stochastic quadratic function and an $\ell_2$-regularized logistic regression problem are presented.
In all the cases tested, our approach improves the convergence of stochastic gradient descent, compensating for the overhead of solving the log posterior maximization.
In particular, pre-conditioning the stochastic gradient with the inverse of our Hessian approximation becomes more advantageous the larger the condition number of the problem is.

\noindent\textbf{AMS subject classifications:}
$\cdot$
65K10 
$\cdot$
90C15 
$\cdot$
62F15 
$\cdot$
90C53 
$\cdot$

\end{abstract}

\maketitle



\section{Introduction}
%
%


First-order optimization methods are a common choice for performing local search when gradients are available.
In the strongly convex case, the gradient descent method converges linearly, requiring a number of iterations proportional to the condition number \cite[Theorem 2.1.15]{nesterov2013introductory}, that is, the ratio between the upper and lower bounds of the eigenvalues of the Hessian matrix in the search domain.
However, first-order methods may not be a practical choice when the condition number for a problem is large, advancing little towards the optimum each iteration.
One way to mitigate this dependency is to use quasi-Newton methods, which use the secant equation to obtain an approximation of the Hessian matrix of the objective function from subsequent gradient evaluations.
Then, one may precondition the gradient with the inverse of the Hessian approximation, furnishing a Newton-like method without the cost of computing the true Hessian or using a finite differences type of approximation \cite{dennis1977quasi}.


In the case of stochastic optimization, the stochastic gradient descent method (SGD) is the canonical first-order method.
The asymptotic convergence of SGD is dominated by noise in gradient evaluations (cf. \cite[Theorem 3]{polyak1987introduction}), thus controlling the condition number for the problem does not improve convergence.
However, controlling the relative statistical error of gradient estimates using a uniform bound recovers linear convergence per iteration for the minimization of strongly convex functions \cite{polyak1987introduction, MICE}.
In this case, the convergence rate depends on the condition number of the problem; hence preconditioning the gradient with the Hessian inverse might be desirable.
The natural approach is to use quasi-Newton methods to approximate the Hessian matrix, but devising quasi-Newton methods for stochastic optimization is challenging due to its noisy gradient observations. 
Naively fitting the secant equations on noisy gradient observations often leads to catastrophic errors in Hessian approximations \cite[Section 4.3.3]{polyak1987introduction}.

There are many works addressing the development of quasi-Newton methods in stochastic optimization.
One of the first attempts is that of Schraudolph, Yu, and G\"unther \cite{schraudolph2007stochastic} where they adapt the L-BFGS method to the stochastic setting by adding damping parameters that control the effect of noise, furnishing the oL-BFGS method.
Their approach is empirical, and not much can be said about the quality of the Hessian approximation.
Bordes and Pierre \cite{bordes2009sgd} propose a diagonal rescaling matrix for linear support vector machines based on oL-BFGS, providing an analysis of their approach.
In their numerical tests, their method matched oL-BFGS in convergence per epoch but with a significant improvement in terms of runtimes.
Hennig \cite{hennig2013fast} uses a Bayesian approach to learn about an objective function from noisy gradient evaluations.
Byrd et al.~\cite{byrd2011use} use a sample average approximation to compute the search direction in Newton-CG and L-BFGS.
As a consequence, they do not have a statistical error on their estimates but a bias.
Hennig \cite{hennig2013fast} models the prior distribution of the objective function as a Gaussian process with a squared exponential kernel and, using the linearity of Gaussian processes, derives prior distributions for the gradient and the Hessian of the objective function.
Then, assuming Gaussian likelihoods for both the gradient and the Hessian of the objective function, Hennig \cite{hennig2013fast} deduces closed-form equations for the mean of the (also Gaussian) posterior distribution.
Sohl-Dickstein, Poole, and Ganguli \cite{sohl2014fast} develop a quasi-Newton method for sums of functions where they keep a quadratic approximation for each of these functions, subsampling them at each iteration and updating the respective quadratic approximation with new gradient information.
They use a shared, low-dimensional subspace to avoid the cost of storing a Hessian matrix for each function.
Byrd et al.\cite{byrd2016stochastic} propose a stochastic quasi-Newton method comprising a two-loop scheme where the outer loop updates the curvature pairs in an L-BFGS scheme (with large sample sizes), and the inner loop performs preconditioned stochastic gradient descent (with lower sample sizes).
Thus, the method proposed by Byrd et al. \cite{byrd2016stochastic} does not incorporate the information from the noisy gradients of the inner loops in the Hessian approximation.
Moritz, Nishihara, and Jordan \cite{moritz2016linearly} improve the stochastic quasi-Newton method proposed by Byrd et al. \cite{byrd2016stochastic} by combining it with a control variate method to estimate the gradient at each iteration, achieving linear convergence.

Gower, Goldfarb, and Richt\'arik \cite{gower2016stochastic} introduce a block BFGS method to update the Hessian approximation each iteration using a sample independent of the sample used to compute the gradient.
The proposed block update fits the Hessian inverse to a sketched Hessian computed from a sub-sample of the data, but requires a Cholesky decomposition each iteration.
Wang et al. \cite{wang2017stochastic} propose a damped quasi-Newton method for stochastic optimization that uses all noisy gradient evaluations to approximate the Hessian inverse.
They provide a convergence proof for nonconvex problems, provided the Hessian eigenvalues have uniform lower bound.
Li \cite{li2017preconditioned} preconditions the stochastic gradient estimates using a matrix computed to fit the secant equation while avoiding amplifying gradient estimate noise.
Bollapragada et al.\cite{bollapragada2018progressive} use standard L-BFGS equations in stochastic optimization by increasing sample sizes to keep the relative error for the gradients uniformly bounded.
Bollapragada et al.\cite{bollapragada2018progressive} recommend a stochastic line-search procedure to be used with their method.
Goldfarb, Ren, and Bahamou~\cite{goldfarb2020practical} develop a quasi-Newton method for deep neural networks that approximates Hessian inverses as block matrices such that each block is decomposed as the kronecker products of smaller matrices.
Lately, Shi et al.~\cite{shi2022noise} build on top of their analysis in \cite{xie2019analysis} to build a quasi-Newton method for stochastic optimization with a noise-resilient Wolfe--Armijo line search.
Their method, however, requires the gradient estimate to be bounded above in norm by a known constant, which they use in their line search.
We refer the interested reader to \cite{bottou2018optimization, mokhtari2020stochastic} for extensive reviews of quasi-Newton methods for stochastic optimization.


Most approaches to tackle the condition number issue in stochastic optimization use equations from deterministic quasi-Newton methods and somehow circumvent the randomness.
Here, we propose a Bayesian framework to build a posterior distribution for the Hessian matrix from noisy gradient observations.
In Bayesian inference, a posterior distribution is built from prior and likelihood distributions.
We model the prior distribution for the Hessian matrix with a 
probability density function (pdf)
that decreases exponentially as the Hessian moves farther away (in the Frobenius norm sense) from the initial Hessian estimate.
This approach is motivated by classic deterministic quasi-Newton methods, where a Frobenius norm regularizer centered at the current Hessian is used to circumvent ill-posedness \cite{dennis1977quasi}.
Moreover, we impose constraints on the Hessian matrix eigenvalues extremes, enforcing positive-definiteness and gradient smoothness.
We model the likelihood using the secant equation, i.e., how likely it is to obtain the observed gradient differences for a given Hessian.
The Hessian that maximizes the posterior distribution, i.e., the Maximum a Posteriori (MAP), is defined as the Hessian approximation.

The choice of prior and likelihood distributions guarantees that the negative log posterior is strongly-convex.
Therefore, we opt to use the Newton-CG method to find the MAP.
Although it is possible to compute the second-order derivatives of the negative log posterior of the Hessian directly, it could be expensive and memory intensive.
Keep in mind that, for a $d$-dimensional problem, the Hessian matrix is a $d \times d$ matrix, thus, the second-order derivatives of the log posterior must be gathered in a fourth-order tensor with $d^4$ components.
Instead of assembling this tensor, we compute its action over any symmetric matrix of adequate dimensions and then use the conjugate gradient method to find the Newton direction.
Moreover, to enforce the stability of Newton-CG, we use the central-path interior point method to impose the prior conditions on the eigenvalues of the Hessian approximation, i.e., we start with a large penalization parameter and decrease it in steps, keeping the problem well-conditioned.
If properly tuned, our approach not only guarantees that Newton-CG converges to the MAP as it also results in a slow decay for the Hessian approximation's smallest eigenvalue, thus improving the stability for the SGD preconditioned with our Hessian approximation.

Other works address the problem of approximating a Hessian from noisy gradient observations.
Hennig \cite{hennig2013fast} also uses Bayesian inference; however, to build a posterior distribution of the objective function as a Gaussian process.
The main difference between the proposed approach and that of Hennig \cite{hennig2013fast} is in the way he models the prior and likelihood distributions.
Moreover, Hennig \cite{hennig2013fast} assumes independence between the components of the gradient noise in his analysis,  and also that they have the same variance.
In the same spirit of the present work, Li \cite{li2017preconditioned} also solves an optimization sub-problem each iteration to find a Hessian approximation from noisy gradients.
His approach differs from ours in the way he uses the available information to estimate the preconditioning matrix; here, we use a Bayesian approach to build a posterior distribution of the Hessian matrix while \cite{li2017preconditioned} obtains a point estimate by minimizing a different objective function.

To illustrate the effectiveness of our approach, we solve two stochastic optimization problems.
The first is minimizing a quadratic function, and the second is training an $\ell_2$ regularized logistic regression model.
In all cases, we start with SGD preconditioned with a diagonal matrix and update the preconditioning matrix in intervals using the proposed Bayesian approach.
We study cases with up to $300$ design variables and condition numbers of orders up to $10^9$.
We present results combining the proposed Bayesian Hessian approach with different variations of SGD, namely, vanilla SGD, SGD-MICE \cite{MICE}, SVRG \cite{johnson2013accelerating}, and SARAH \cite{nguyen2017sarah}.

\subsection{Problem definition}

The stochastic optimization problem is stated as follows. Let $\bs{\xi}$ be the design variable in dimension $d_{\bs{\xi}}$ and $\bs{\theta}$ be a vector-valued random variable in dimension $d_{\bs{\theta}}$ equipped with a probability density function $\varrho\colon\rset^{d_{\bs{\theta}}}\to \rset^+$.
Here, $\bb{E}[\cdot|\bs{\xi}]$ and $\bb{V}[\cdot|\bs{\xi}]$ are the expectation and variance operators conditioned on $\bs{\xi}$, respectively. Aiming at minimizing the expectations conditioned on $\bs{\xi}$, we state the problem as
\begin{equation}\label{eq:underlying.problem}
\bs{\xi}^* = \underset{\bs{\xi} \in \Xi}{\arg\min} \, \bb{E} [f(\bs{\xi},\bs{\theta})|\bs{\xi}],
\end{equation}
where $\Xi\subset\bb{R}^{d_{\bs{\xi}}}$ is a feasible set and $f \colon \Xi \times \rset^{d_{\bs{\theta}}} \to \bb{R}$.
Through what follows, let $F(\bs{\xi})\coloneqq\bb{E} [f(\bs{\xi},\bs{\theta})|\bs{\xi}]$ where $F \colon \Xi \to \bb{R}$.

In minimizing \eqref{eq:underlying.problem} with respect to the design variable $\bs{\xi} \in \Xi$, vanilla SGD is constructed with the update rule
\begin{equation}\label{eq:sgd}
\bs{\xi}_{k+1} = \bs{\xi}_{k} - \eta_k \bs{\upsilon}_k,
\end{equation}
where $\eta_k$ is the step size at iteration $k$, and $\bs{\upsilon}_k$ is an unbiased estimator of the gradient of $F$ at $\bs{\xi}_k$.
For instance, an unbiased estimator $\bs{\upsilon}_k$ of the gradient of $F$ at $\bs{\xi}_k$ may be constructed using a Monte Carlo estimator,
\begin{equation}\label{eq:grad.def.2}
\nabla_{\bs{\xi}} F(\bs{\xi}_k) =  \bb{E}[f(\bs{\xi}_k,\bs{\theta})|\bs{\xi}_k] \approx \bs{\upsilon}_k :=\dfrac{1}{b} \sum_{i=1}^b \nabla_{\bs{\xi}} f(\bs{\xi}_k,\bs{\theta}_i),
\end{equation}
with $b$ independent and identically distributed (iid) random variables $\bs{\theta}_i\sim\varrho$.
The estimator \eqref{eq:grad.def.2} is a random variable and its use in optimization algorithms gives rise to the so-called \emph{Stochastic Optimizers}.

First-order methods are the most popular choice in stochastic optimization due to their robustness to noise.
However, these methods might suffer from slow convergence when the condition number of the problem is large.
Here, we are interested in finding a matrix $\hat{\bs{B}}_k$ that approximates $\nabla^2 F(\bs{\xi}_k)$ so that its inverse can be used as a preconditioner for the stochastic gradient,
\begin{equation}
  \bs{\xi}_{k+1} = \bs{\xi}_{k} - \eta_k \hat{\bs{B}}_k^{-1} \bs{\upsilon}_k,
\end{equation}
thus diminishing the effect of the condition number on the optimization convergence.
In the deterministic setting, quasi-Newton methods circumvent the sensibility to the condition number by using the gradient observations to construct the Hessian matrix approximation at the current iterate.
However, in the stochastic optimization setting, we first have to address the following questions:
\begin{enumerate}
  \item[\textbf{Q1}:] How does one obtain a Hessian approximation from noisy gradient observations?
  \item[\textbf{Q2}:] How do we avoid amplifying the noise in gradient estimates with our preconditioner matrix $\hat{\bs{B}}_k^{-1}$? 
  \item[\textbf{Q3}:] To mitigate the effect of the condition number, is it enough to have an accurate Hessian approximation, however, acting on noisy gradient estimates?
\end{enumerate}

To address \textbf{Q1}, consider the curvature pair $\bs{y}_k \coloneqq \nabla_{\bs{\xi}}F(\bs{\xi}_k) - \nabla_{\bs{\xi}}F(\bs{\xi}_{k-1})$ and $\bs{s}_k \coloneqq \bs{\xi}_k - \bs{\xi}_{k-1}$. Thus, we first aim to find the Hessian approximation that minimizes the residue of the secant equation, that is,
\begin{equation} \label{eq:problem}
\min_{\bs{B} \in \bb{R}^{d_{\bs{\xi}}}\times\bb{R}^{d_{\bs{\xi}}}} \sum_{\ell=1}^k \| \bs{B} \bs{s}_\ell - \bs{y}_\ell \|^2.
\end{equation}
However, we do not have access to $\bs{y}_k\coloneqq\bb{E}[\nabla_{\bs{\xi}}f(\bs{\xi}_k, \bs{\theta})|\bs{\xi}_k] - \bb{E}[\nabla_{\bs{\xi}}f(\bs{\xi}_{k-1}, \bs{\theta})|\bs{\xi}_{k-1}]$,
instead, we may sample finitely many times $\hat{\bs{y}}_k(\bs{\theta}) \coloneqq \nabla_{\bs{\xi}}f(\bs{\xi}_k, \bs{\theta}) - \nabla_{\bs{\xi}}f(\bs{\xi}_{k-1}, \bs{\theta})$.
In Section \ref{sec:bayes}, we discuss how to use Bayesian inference to find a matrix $\hat{\bs{B}}_k$ given the noisy curvature pairs $\hat{\bs{y}}$ and $\bs{s}$.

One of the main problems faced by preconditioning in stochastic optimization is the amplification of the noise in gradient estimates.
This amplification is proportional to the smallest eigenvalue of matrix $\hat{\bs{B}}_k$,
\begin{equation} \label{eq:noise}
\bb{E}\left[
  \left\|
    \hat{\bs{B}}_k^{-1} \bs{\upsilon}_k
   - \bb{E}\left[
      \hat{\bs{B}}_k^{-1} \bs{\upsilon}_k
      \Big|\hat{\bs{B}}_k^{-1}
    \right]
  \right\|^2
  \bigg| \hat{\bs{B}}_k^{-1}
\right] \le
\left\|
  \hat{\bs{B}}_k^{-1}
\right\|^2
\bb{E}\left[
  \left\|\bs{\upsilon}_k - \bb{E}[\bs{\upsilon}_k]
  \right\|^2
  \Big|\hat{\bs{B}}_k^{-1}
\right],
\end{equation}
by the matrix-norm compatibility.
Thus, to tackle \textbf{Q2}, we impose a lower bound constraint on the eigenvalues of $\hat{\bs{B}}_k$, i.e., $\hat{\bs{B}}_k \succeq \tilde{\mu}\bs{I}$.

To address \textbf{Q3} in a numerical fashion, we investigate and compare stochastic gradients preconditioned with our Hessian inverse approximation with and without control in the relative statistical error.

Regarding notation, in this work, we use bold-face lower-case letters as vectors and bold-face upper-case letters as matrices.
Moreover, we refer to first-order derivatives of functions with respect to matrices as gradients, even though they are also matrices.



\section{Approximating the Hessian and its inverse using Bayesian inference}
\label{sec:bayes}

Important to what follows is the Monte Carlo estimator of $\bs{y}_k$
\begin{equation} \label{eq:y_bar}
  \bar{\bs{y}}_k \coloneqq \frac{1}{b_k} \sum_{i=1}^{b_k} \hat{\bs{y}}_k(\bs{\theta}_i),
\end{equation}
with $b_k$ iid random variables $\bs{\theta}_i\sim\varrho$.

\subsection{Bayesian formulation of the Hessian inference problem}
For the pair of sets $Y \coloneqq \{\bar{\bs{y}}_\ell\}_{\ell=1}^k$ and $S \coloneqq \{\bs{s}_\ell\}_{\ell=1}^k$, we employ the Bayes' formula to build a posterior distribution of the Hessian $\bs{B}$ given $Y$ and $S$,
\begin{equation} \label{eq:bayes}
  \pi(\bs{B} | S, Y) \propto p(S, Y | \bs{B}) \pi(\bs{B}),
\end{equation}
where the likelihood distribution $p(S, Y | \bs{B})$ and the prior distribution $\pi(\bs{B})$ will be discussed in sections \ref{sec:likelihood} and \ref{sec:prior}, respectively.
We opt to use as a Hessian approximation the matrix that minimizes the negative log posterior subject to the eigenvalues constraints
\begin{equation}\label{eq:subproblem}
\bs{B}_k^* \coloneqq
\left\{
\begin{aligned}
& \arg\min_{\bs{B}}
\left(
\cl{L}_k(\bs{B}) \coloneqq -\log(\pi(\bs{B} | S, Y))
\right) \\[4pt]
& \text{subject to }
\begin{aligned}[t]
\lambda_{\min}(\bs{B}) &\ge \tilde{\mu}\\
\lambda_{\max}(\bs{B}) &\le L\\
\bs{B} &= \bs{B}^{\trans},
\end{aligned}
\end{aligned}
\right.
\end{equation}
where $\lambda_{min}$ and $\lambda_{max}$ are operators that return the minimum and maximum eigenvalues, respectively; $L$ is the smoothness constant of $F$; and $\tilde{\mu} > 0$ is an arbitrary parameter.
If $F$ is $\mu$-convex, $\tilde{\mu}$ can be chosen as $\mu$, however, a larger $\tilde{\mu}$ might be desirable to avoid amplifying the noise in the search direction.
Since computing $\bs{B}^*_k$ is not possible in general, we use an approximation $\hat{\bs{B}}_k$, obtained numerically, as a Hessian approximation.

\subsubsection{
Likelihood of observing curvature pairs given a Hessian 
}
\label{sec:likelihood}

We use the secant equations to derive the 
likelihood of observing the curvature pairs $(S, Y)$ given a Hessian $\bs{B}$.
Let the difference in gradients be as in \eqref{eq:y_bar}.
From the observations $\hat{\bs{y}}_k$, we build a sample covariance matrix $\bs{\Sigma}_k$.
Then, the covariance matrix of the estimator $\bar{\bs{y}}_k$ is given by
\begin{equation}
  \bar{\bs{\Sigma}}_k \coloneqq 
  \textrm{Cov}
  [\bar{\bs{y}}_k] = \frac{\bs{\Sigma}_k}{b_k}.
\end{equation}
We are interested in the inverse of $\bar{\bs{\Sigma}}_k$.
However, to avoid singularity issues, we use a precision matrix approximation $\bs{P}_{\ell}$ instead of $\bar{\bs{\Sigma}}_\ell^{-1}$
for each $\bar{\bs{y}}_\ell \in Y$.
We define $\bs{P}_{\ell}$ as
\begin{equation}
  \bs{P}_\ell \coloneqq
  \left(
    \bar{\bs{\Sigma}}_\ell + \sigma_P (\max(\bar{\bs{\Sigma}}_\ell)) \bs{I}_{d_{\xi}}
  \right) ^{-1},
\end{equation}
where $\sigma_P$ is a regularization parameter and $\bs{I}_{d_{\xi}}$ is the identity matrix of size $d_{\xi}$.
Then, we define the negative log likelihood of curvature pair $\bs{s}_k$ and $\bar{\bs{y}}_k$ given a matrix $\bs{B}$ as
\begin{equation}
  - \log(p(\bs{s}_k, \bar{\bs{y}}_k | \bs{B})) = \frac{1}{2 \nu} \| \bs{B} \bs{s}_k - \bar{\bs{y}}_k \|_{\bs{P}_k}^2 + C_k,
\end{equation}
thus, for $S$ and $Y$ generated by a first-order stochastic gradient descent method,
\begin{equation} \label{eq:log_lkl}
  - \log(p(S, Y | \bs{B})) = \frac{1}{2 \nu} \sum_{\ell=1}^k \| \bs{B} \bs{s}_\ell - \bar{\bs{y}}_\ell \|_{\bs{P}_\ell}^2 + C.
\end{equation}
In the above equations, 
$C = \sum_{i=\ell}^k C_\ell$
are constants and
$\nu$ is used to normalize the negative log likelihood,
\begin{align}
  \label{eq:beta_hat}
  \nu = \sum_{\ell=1}^k \|\bs{P}_\ell (\hat{\bs{B}}_{k-1} \bs{s}_\ell - \bar{\bs{y}}_\ell) \|
  \| \bs{s}_\ell \|,
\end{align}
being $\hat{\bs{B}}_{k-1}$ a Hessian approximation at iteration $k-1$.

\subsubsection{Prior distribution of the Hessian matrix} \label{sec:prior}

The prior distribution must encode the knowledge available about the Hessian before the observations of $S$ and $Y$.
In spirit of quasi-Newton methods, we want to discourage a Hessian that differs much from the current estimate.
For this reason, we adopt a negative log prior of $\bs{B}$ as the squared weighted Frobenius norm of the difference between $\bs{B}$ and $\hat{\bs{B}}_{k-1}$,
\begin{equation} \label{eq:frob_reg}
  \frac{1}{2} \|\bs{B} - \hat{\bs{B}}_{k-1} \|_{F, \bs{W}}^2 =
  \frac{1}{2} \|\bs{W}^{\nicefrac{1}{2}} (\bs{B} - \hat{\bs{B}}_{k-1}) \bs{W}^{\nicefrac{1}{2}} \|_F^2,
\end{equation}
for a nonsingular symmetric matrix of weights $\bs{W}$.
Another information available about the matrix $\bs{B}$ is that its eigenvalues must lay between $\tilde{\mu}$ and $L$, as described in the constraints of \eqref{eq:subproblem}.
To enforce these eigenvalue constraints on $\bs{B}$, we define the prior distribution from its negative logarithm as
\begin{equation}
  - \log (\pi(\bs{B})) =
  \begin{cases}
    \frac{1}{2} \|\bs{B} - \hat{\bs{B}}_{k-1} \|_{F, \bs{W}}^2 + C  & \text{if }
    \hat{\mu} \bs{I}_{d_{\xi}} \preceq \bs{B} \preceq \hat{L} \bs{I}_{d_{\xi}}
    \\
    + \infty                                          & \text{otherwise},
  \end{cases}
\end{equation}
where again $C$ is a normalizing constant.
Here, we opt to impose the eigenvalues constraints on the negative log prior as logarithmic barriers of form
\begin{align}
  - \log\left(\det\left(\bs{B} - \hat{\mu} \bs{I}_{d_{\xi}}\right)\right)
  - \log\left(\det\left(\hat{L} \bs{I}_{d_{\xi}} - \bs{B}\right)\right),
\end{align}
where
\begin{equation}
  \hat{\mu} \coloneqq \frac{\tilde{\mu}}{\alpha},
  \quad
  \hat{L} \coloneqq L \alpha,
\end{equation}
and $\alpha \geq 1$ is a relaxation parameter.
Hence, we propose a family of priors with logarithmic barriers with hyperparameter $\beta$ as
\begin{multline} \label{eq:log_prior}
  - \log(\pi(\bs{B};\beta)) =  \frac{\rho}{2} \|\bs{B} - \hat{\bs{B}}_{k-1} \|_{F, \bs{W}}^2
  - \beta \log\left(\det\left(\bs{B} - \hat{\mu} \bs{I}_{d_{\xi}}\right)\right) \\
  - \beta \log\left(\det\left(\hat{L} \bs{I}_{d_{\xi}} - \bs{B}\right)\right) + C,
\end{multline}
being supported on the set of symmetric positive-definite matrices with eigenvalues between $\hat{\mu}$ and $\hat{L}$.

The parameter $\rho > 0$ weights the Frobenius regularization term and $\beta > 0$ weights the logarithmic barrier terms.
Note that $\log (\pi(\bs{B}; \beta)) \rightarrow \log(\pi(\bs{B}))$ as $\beta \downarrow 0$ with $\rho=1$.

Some options for $\bs{W}$ are $\bs{I}_{d_{\xi}}$ and $\hat{\bs{B}}_{k-1}^{-\nicefrac{1}{2}}$.
According to \cite{dennis1977quasi}, the Frobenius norm as in \eqref{eq:frob_reg} with weights $\hat{\bs{B}}_{k-1}^{-\nicefrac{1}{2}}$ measures the relative error of $\bs{B}$ with respect to $\hat{\bs{B}}_{k-1}$.
In practice, estimating $\hat{\bs{B}}_{k-1}^{-\nicefrac{1}{2}}$ can be expensive.
In these cases, we recommend using $\bs{W} = \bs{I}_{d_{\xi}}$.

\subsubsection{Posterior distribution of the Hessian matrix}

Substituting \eqref{eq:log_lkl} and \eqref{eq:log_prior} in \eqref{eq:bayes}, we obtain the posterior distribution of the Hessian matrix $\bs{B}$ given curvature pairs $S$ and $Y$.
The negative log posterior with parameter $\beta$ is
\begin{multline} \label{eq:log_post}
  \cl{L}_k(\bs{B}; \beta) \coloneqq \frac{1}{2 \nu} \sum_{\ell=1}^k \| \bs{B} \bs{s}_\ell - \bar{\bs{y}}_\ell \|_{\bs{P}_\ell}^2
  + \frac{\rho}{2} \|\bs{B} - \hat{\bs{B}}_{k-1} \|_{F, \bs{W}}^2
  - \beta \log\left(\det\left(\bs{B} - \hat{\mu} \bs{I}_{d_{\xi}}\right)\right) \\
  - \beta \log\left(\det\left(\hat{L} \bs{I}_{d_{\xi}} - \bs{B}\right)\right) + C.
\end{multline}
The equation \eqref{eq:log_post} can be interpreted as the sum of the secant equation residues for each curvature pair (in the metrics induced by $\bs{P}_\ell$), subject to the logarithmic barrier terms that impose the eigenvalues constraints, and a Frobenius norm regularization term.

Since we are interested in minimizing the negative log posterior in \eqref{eq:log_post}, it is our interest to know if it is strongly convex.
The negative log-likelihood is convex but not strictly-convex until enough linearly independent curvature pairs  $(\bs{s}_\ell, \bs{y}_\ell)$ are observed.
Regarding the logarithmic barrier terms of the prior, the negative logarithm of the determinant defined on symmetric positive-definite matrices is a convex function \cite[Section 3.1.5]{boyd2004convex}; however, not a strictly-convex one.
The Frobenius norm regularization term with $\bs{W} = \bs{I}_{d_{\xi}}$ is $\rho$-convex, therefore the negative log posterior is strongly convex with parameter $\rho$.
Here, the regularization term works as a Tikhonov regularization term ($\|\bs{A}\|_F \equiv \|\text{vec}(\bs{A})\|_2$), guaranteeing uniqueness of solution.
Properly tuning $\beta$ and $\rho$ guarantees that the problem of minimizing \eqref{eq:log_post} is well-conditioned.
A discussion on how to choose $\beta$ is presented in Section \ref{sec:find_map}.

\subsection{Finding the posterior maximizer} \label{sec:find_map}

Finding the $\bs{B}^*_k$ that solves the problem in \eqref{eq:subproblem} is not a trivial task.
Since the negative log posterior is $\rho$-convex, local search methods like the Newton method can be used to find its minimizer.
However, the Newton method requires assembling the fourth-order tensor of second-order derivatives of the negative log posterior, which can be expensive in terms of processing and memory allocation.
Instead, we opt to derive the action of this fourth-order tensor of second-order derivatives over any given matrix and use the conjugate gradient method to find the Newton direction, a procedure known as the Newton-CG method.
Since we have a constraint in \eqref{eq:subproblem} that $\bs{B}$ should be symmetric, the first-order optimality condition is
\begin{equation}
  \sym( \nabla_{\bs{B}} \cl{L}_k(\bs{B}; \beta)) = 0,
\end{equation}
where the $\sym(\cdot)$ operator returns the symmetric part of its argument, e.g., $\sym(\bs{A}) = \frac{1}{2} (A + A^T)$.
Thus, to use the Newton-CG method to solve this problem, every iteration we need to find the direction $\bs{D}$ that satisfies
\begin{equation}
  \sym (\delta_{\bs{D}} \nabla_{\bs{B}} \cl{L}_k(\bs{B}; \beta))
  = - \sym(\nabla_{\bs{B}} \cl{L}_k(\bs{B}; \beta)),
\end{equation}
where $\delta_{\bs{D}}$ is the directional derivative with respect to $\bs{B}$ in the direction of $\bs{D}$,
\begin{equation}
  \delta_{\bs{D}} g(\bs{B}) =
  \lim_{\vartheta \downarrow 0} \frac{
      g(\bs{B}+\vartheta \bs{D}) - g(\bs{B})
    }{
      \vartheta
    }.
\end{equation}

To avoid ill-conditioning, we need a large enough regularization parameter $\rho$, but not too large since it introduces bias in the posterior maximizer.
Also, the bias given by the log-barriers grows with $\beta$, the logarithmic barrier parameter; but a small $\beta$ negatively impacts the smoothness of the posterior distribution.
We opt to keep $\rho$ fixed and use an interior-point central-path method \cite[Section 11.2.2]{boyd2004convex} to decrease the parameter $\beta$; the solution $\hat{\bs{B}}_k(\beta)$ of the approximated problem
\begin{equation}
  \hat{\bs{B}}_k(\beta) = \arg\min_{\bs{B}}
  \left(
  \cl{L}_k(\bs{B}; \beta) \right)
\end{equation}
converges to the solution $\bs{B}_k^*$ of \eqref{eq:subproblem} as $\beta \downarrow 0$.
Thus, we start with a large $\beta$ and a large tolerance $tol$ and perform optimization until the norm of the gradient of the log posterior is below $tol$. Then, we decrease both $\beta$ and $tol$ by a factor $\gamma$ and repeat the process.
This procedure is repeated until the desired tolerance is achieved with the desired $\beta$.
To minimize the negative log posterior for a given $\beta$, we use the Newton-CG method with backtracking-Armijo line search; we use the conjugate gradient method to find the Newton direction and then the backtracking-Armijo to find a step size.
Therefore, to use the Newton-CG, we need the gradient of $\cl{L}_k$ with respect to $\bs{B}$ and the directional derivative of the gradient of $\cl{L}_k$ with respect to $\bs{B}$ in the direction of any arbitrary symmetric matrix $\bs{V} \in \bb{R}^{d_{\bs{\xi}}} \times \bb{R}^{d_{\bs{\xi}}}$.
The pseudo-code for the minimization of the negative log posterior is presented in Algorithm \ref{alg:find_map}.
We use the central-path method to decrease $\beta$ and $tol$ every step outer loop, indexed as $i$ in Algorithm 1.
In lines 6 to 10, we have the Newton-CG method, where, each iteration, the conjugate gradient method is used to find the approximate Newton direction $\bs{D}$ and a backtracking-Armijo line-search is used to find the step size.

\begin{algorithm}[h]
 \setstretch{1.5}
 \caption{Pseudocode of the procedure to find the maximum a posteriori of the Hessian matrix of the objective function given noisy curvature pair observations.
 }\label{alg:find_map}
 \begin{algorithmic}[1]
  \Procedure{FindMAP}{$S, Y, \hat{\bs{B}}_{k-1}, \beta_0, tol_0, \gamma=2, c_A=10^{-4}$}
  \State $\hat{\bs{B}}_k \gets \hat{\bs{B}}_{k-1}$
  \For{$i=1, 2,...$}
  \State $\beta_i = \beta_{i-1}/\gamma$
  \State $tol_i = tol_{i-1}/\gamma$
  \While{$\|\nabla_{\bs{B}} \cl{L}_k(\hat{\bs{B}}_k; \beta_i)\|_F > tol_i$}
  \State use the conjugate gradient method to find a symmetric $\bs{D}$ such that \\
  \begin{flushright}
  $\sym (\delta_{\bs{D}} \nabla_{\bs{B}} \cl{L}_k(\hat{\bs{B}}_k; \beta_i)) = -\sym(\nabla_{\bs{B}} \cl{L}_k(\hat{\bs{B}}_k; \beta_i))$
  \end{flushright}
  \State find $\varsigma$ such that $\cl{L}_k(\hat{\bs{B}}_k + \varsigma \bs{D}; \beta_i) \le \cl{L}_k(\hat{\bs{B}}_k; \beta_i) + c_A \varsigma \nabla_{\bs{B}} \cl{L}_k(\hat{\bs{B}}_k; \beta_i) \colon \bs{D}$ using backtracking-Armijo
  \State $\hat{\bs{B}}_k \gets \hat{\bs{B}}_k + \varsigma \bs{D}$
  \EndWhile
  \EndFor
  \State \textbf{return} $\hat{\bs{B}}_k$
  \EndProcedure
 \end{algorithmic}
\end{algorithm}

Next, we will present how to compute $\nabla_{\bs{B}} \cl{L}_k(\hat{\bs{B}}_k; \beta)$ and $\delta_{\bs{V}} \nabla_{\bs{B}} \cl{L}_k(\hat{\bs{B}}_k; \beta)$ for any symmetric $\bs{V} \in \bb{R}^{d_{\bs{\xi}}} \times \bb{R}^{d_{\bs{\xi}}}$.
Note that, since $\bs{B}$ is a matrix, $\nabla_{\bs{B}} \cl{L}_k$ and $\delta_{\bs{V}} \nabla_{\bs{B}} \cl{L}_k$ are also symmetric matrices of the same dimensions of $\bs{B}$.

\subsection{Gradient of the negative log posterior}
To compute $\nabla_{\bs{B}} \cl{L}_k$, we need to compute the gradients of the log prior and of the log likelihood,
\begin{equation} \label{eq:grad_log_post_}
  \nabla_{\bs{B}} \cl{L}_k(\bs{B}; \beta) = -\nabla_{\bs{B}} \log (p(S, Y | \bs{B})) -\nabla_{\bs{B}} \log (\pi(\bs{B}; \beta)).
\end{equation}

The gradient of the negative log likelihood can be calculated from \eqref{eq:log_lkl} as
\begin{align}
  - \nabla_{\bs{B}} \log (p(S, Y | \bs{B})) &= \frac{1}{2 \nu} \sum_{\ell=1}^k \nabla_{\bs{B}} \| \bs{B} \bs{s}_\ell - \bar{\bs{y}}_\ell \|_{\bs{P}_\ell}^2 \\
  &= \label{eq:grad_lkl}
   \frac{1}{\nu} \sum_{\ell=1}^k \bs{P}_\ell (\bs{B} \bs{s}_\ell - \bar{\bs{y}}_\ell) \otimes \bs{s}_\ell.
\end{align}

The gradient of the negative log prior can be derived from \eqref{eq:log_prior} as
\begin{multline} \label{eq:grad_prior_}
  -\nabla_{\bs{B}} \log (\pi(\bs{B}; \beta))
  = \frac{\rho}{2} \nabla_{\bs{B}} \|\bs{W}^{\nicefrac{1}{2}} (\bs{B} - \hat{\bs{B}}_{k-1}) \bs{W}^{\nicefrac{1}{2}} \|_F^2
  - \beta \nabla_{\bs{B}} \log\left(\det\left(\bs{B} - \hat{\mu} \bs{I}_{d_{\xi}}\right)\right) \\
  - \beta \nabla_{\bs{B}} \log\left(\det\left(\hat{L} \bs{I}_{d_{\xi}} - \bs{B}\right)\right).
\end{multline}
The gradient of the Frobenius norm regularizer is
\begin{align} \label{eq:grad_prior_1}
  \frac{1}{2} \nabla_{\bs{B}} \|\bs{W}^{\nicefrac{1}{2}} (\bs{B} - \hat{\bs{B}}_{k-1}) \bs{W}^{\nicefrac{1}{2}} \|_F^2 &=
  \bs{W} (\bs{B} - \hat{\bs{B}}_{k-1}) \bs{W},
\end{align}
whereas the gradients of the logarithmic barrier terms are derived using the Jacobi formula,
\begin{align}
  \label{eq:grad_prior_2a}
  - \nabla_{\bs{B}} \log\left(\det\left(\bs{B} - \hat{\mu} \bs{I}_{d_{\xi}}\right)\right)
  &= - (\bs{B} - \hat{\mu}\bs{I}_{d_{\xi}})^{-1} \\
  \label{eq:grad_prior_2b}
  - \nabla_{\bs{B}} \log\left(\det\left(\hat{L} \bs{I}_{d_{\xi}} - \bs{B}\right)\right)
  &=  (\hat{L}\bs{I}_{d_{\xi}} - \bs{B})^{-1}.
\end{align}
Substituting \eqref{eq:grad_prior_1}, \eqref{eq:grad_prior_2a}, and \eqref{eq:grad_prior_2b} in \eqref{eq:grad_prior_},
\begin{equation} \label{eq:grad_prior}
-\nabla_{\bs{B}} \log (\pi(\bs{B}; \beta))
= \rho \bs{W} (\bs{B} - \hat{\bs{B}}_{k-1}) \bs{W}
- \beta (\bs{B} - \hat{\mu}\bs{I}_{d_{\xi}})^{-1}
+ \beta (\hat{L}\bs{I}_{d_{\xi}} - \bs{B})^{-1}.
\end{equation}
Finally, substituting the gradient of the negative log likelihood \eqref{eq:grad_lkl} and the gradient of the negative log prior \eqref{eq:grad_prior} in \eqref{eq:grad_log_post_},
we obtain the gradient of the negative log posterior
\begin{multline} \label{eq:grad_log_post}
  \nabla_{\bs{B}} \cl{L}_k(\bs{B}; \beta) =
  \frac{1}{\nu}
  \sum_{\ell=1}^k \bs{P}_\ell (\bs{B} \bs{s}_\ell - \bar{\bs{y}}_\ell) \otimes \bs{s}_\ell
  + \rho \bs{W} (\bs{B} - \hat{\bs{B}}_{k-1}) \bs{W} \\
  - \beta (\bs{B} - \hat{\mu}\bs{I}_{d_{\xi}})^{-1} + \beta (\hat{L}\bs{I}_{d_{\xi}} - \bs{B})^{-1}.
\end{multline}

\subsection{Directional derivative of the gradient of the negative log posterior}
As discussed in Section \ref{sec:find_map}, to find the posterior distribution maximizer using Newton-CG, we need to compute the directional derivative of the gradient of the negative log posterior in a given direction.
Thus, for any symmetric matrix $\bs{V} \in \bb{R}^{d_{\bs{\xi}}} \times \bb{R}^{d_{\bs{\xi}}}$,
we derive $\delta_{\bs{V}} \nabla_{\bs{B}} \cl{L}_k$ as
\begin{equation} \label{eq:dir_hess_post_}
  \delta_{\bs{V}} \nabla_{\bs{B}} \cl{L}_k(\bs{B}; \beta) = -\delta_{\bs{V}} \nabla_{\bs{B}} \log (p(S, Y | \bs{B})) -\delta_{\bs{V}} \nabla_{\bs{B}} \log (\pi(\bs{B}; \beta)).
\end{equation}
Following from \eqref{eq:grad_lkl},
\begin{align}
  -\delta_{\bs{V}} \nabla_{\bs{B}} & \log (p(S, Y | \bs{B})) \nonumber \\
  &= -\lim_{\vartheta \downarrow 0}
    \frac{\nabla_{\bs{B}}  \log p(S, Y | \bs{B} + \vartheta \bs{V}) - \nabla_{\bs{B}} \log (p(S, Y | \bs{B}))}{\vartheta} \\
  &= \frac{1}{\nu} \lim_{\vartheta \downarrow 0}
      \frac{
        \sum_{\ell=1}^k \bs{P}_\ell ((\bs{B}+\vartheta \bs{V}) \bs{s}_\ell - \bar{\bs{y}}_\ell) \otimes \bs{s}_\ell
        - \sum_{\ell=1}^k \bs{P}_\ell (\bs{B} \bs{s}_\ell - \bar{\bs{y}}_\ell) \otimes \bs{s}_\ell
      }{\vartheta} \\
  &= \frac{1}{\nu} \sum_{\ell=1}^k \bs{P}_\ell \bs{V} \bs{s}_\ell \otimes \bs{s}_\ell.
  \label{eq:dir_hess_lkl}
\end{align}
The directional derivative of the gradient of the negative log prior is derived from \eqref{eq:grad_prior},
\begin{multline} \label{eq:dir_hess_pr_}
  -\delta_{\bs{V}} \nabla_{\bs{B}} \log (\pi(\bs{B}; \beta))
  = \rho \, \delta_{\bs{V}}
    \left(\bs{W} (\bs{B} - \hat{\bs{B}}_{k-1}) \bs{W} \right) \\
  - \beta \delta_{\bs{V}} (\bs{B} - \hat{\mu}\bs{I}_{d_{\xi}})^{-1}
  + \beta \delta_{\bs{V}} (\hat{L}\bs{I}_{d_{\xi}} - \bs{B})^{-1}.
\end{multline}
The Frobenius regularizer term from \eqref{eq:dir_hess_pr_} is derived as
\begin{align}
  \delta_{\bs{V}} \left(\bs{W} (\bs{B} - \hat{\bs{B}}_{k-1}) \bs{W} \right)
  &= \lim_{\vartheta \downarrow 0}
    \frac{
    \bs{W} (\bs{B} + \vartheta \bs{V} - \hat{\bs{B}}_{k-1}) \bs{W}
    - \bs{W} (\bs{B} - \hat{\bs{B}}_{k-1}) \bs{W}
    }{
    \vartheta
    }
    \nonumber
   \\
   & = \bs{W} \bs{V} \bs{W}.  \label{eq:dir_hess_pr_1}
\end{align}


Following the same rationale for the log-barrier terms,
\begin{align}
  \delta_{\bs{V}} (\bs{B} &- \hat{\mu}\bs{I}_{d_{\xi}})^{-1} 
  \nonumber \\
  &= \nonumber
  \lim_{\vartheta \downarrow 0}
  \frac{
  ((\bs{B} + \vartheta \bs{V}) - \hat{\mu}\bs{I}_{d_{\xi}})^{-1}
  - (\bs{B} - \hat{\mu}\bs{I}_{d_{\xi}})^{-1}
  }{
  \vartheta
  } \\
  &= \nonumber
  \lim_{\vartheta \downarrow 0}
  \frac{
  ((\bs{B} - \hat{\mu}\bs{I}_{d_{\xi}})
    (\bs{I}_{d_{\xi}} + \vartheta (\bs{B} - \hat{\mu}\bs{I}_{d_{\xi}})^{-1} \bs{V}))^{-1}
  - (\bs{B} - \hat{\mu}\bs{I}_{d_{\xi}})^{-1}
  }{
  \vartheta
  } \\
  &= \nonumber
  \lim_{\vartheta \downarrow 0}
  \frac{
  (\bs{I}_{d_{\xi}} + \vartheta (\bs{B} - \hat{\mu}\bs{I}_{d_{\xi}})^{-1} \bs{V})^{-1}
  (\bs{B} - \hat{\mu}\bs{I}_{d_{\xi}})^{-1}
  - (\bs{B} - \hat{\mu}\bs{I}_{d_{\xi}})^{-1}
  }{
  \vartheta
  } \\  
  &= \nonumber
  \lim_{\vartheta \downarrow 0}
  \frac{
  (\bs{I}_{d_{\xi}} - \vartheta (\bs{B} - \hat{\mu}\bs{I}_{d_{\xi}})^{-1} \bs{V} + \cl{O}(\vartheta^2))
  (\bs{B} - \hat{\mu}\bs{I}_{d_{\xi}})^{-1}
  - (\bs{B} - \hat{\mu}\bs{I}_{d_{\xi}})^{-1}
  }{
  \vartheta
  } \\
  \label{eq:dir_hess_pr_2a}
&= -(\bs{B} - \hat{\mu}\bs{I}_{d_{\xi}})^{-1} 
\bs{V} (\bs{B} - \hat{\mu}\bs{I}_{d_{\xi}})^{-1},
\end{align}
and, 
\begin{align}
  \delta_{\bs{V}} (\hat{L} \bs{I}_{d}-\bs{B})^{-1}
  &= \lim_{\vartheta \downarrow 0}
    \frac{
    ((\hat{L} \bs{I}_{d}-(\bs{B} + \vartheta \bs{V})))^{-1}
    - (\hat{L} \bs{I}_{d}-\bs{B})^{-1}
    }{
    \vartheta
    }
    \nonumber
   \\
   & = \lim_{\vartheta \downarrow 0} 
   \frac{((\hat{L} \bs{I}_{d}-\bs{B}) (\bs{I}_{d_{\xi}}-\vartheta (\hat{L} \bs{I}_{d}-\bs{B})^{-1} \bs{V}))^{-1} - (\hat{L} \bs{I}_{d}-\bs{B})^{-1}}{\vartheta} \nonumber\\[4pt]
   &= \underset{\vartheta \downarrow 0}{\lim}
  \frac{(\bs{I}_{d_{\xi}} - \vartheta (\hat{L} \bs{I}_{d}-\bs{B})^{-1} \bs{V})^{-1} (\hat{L} \bs{I}_{d}-\bs{B})^{-1} - (\hat{L} \bs{I}_{d}-\bs{B})^{-1}}{\vartheta} \nonumber\\[4pt]
   &= \underset{\vartheta \downarrow 0}{\lim}
 \frac{(\bs{I}_{d_{\xi}}+\vartheta (\hat{L} \bs{I}_{d}-\bs{B})^{-1} \bs{V} + \cl{O}(\vartheta^2)) (\hat{L} \bs{I}_{d}-\bs{B})^{-1} - (\hat{L} \bs{I}_{d}-\bs{B})^{-1}}{\vartheta} \nonumber\\[4pt]
   &= 
   \label{eq:dir_hess_pr_2b}
   (\hat{L} \bs{I}_{d}-\bs{B})^{-1} \bs{V} (\hat{L} \bs{I}_{d}-\bs{B})^{-1}.
\end{align}

Substituting \eqref{eq:dir_hess_pr_1}, \eqref{eq:dir_hess_pr_2a}, and \eqref{eq:dir_hess_pr_2b} back into \eqref{eq:dir_hess_pr_}, we get the directional derivative of the gradient of the negative log prior,
\begin{multline} \label{eq:dir_hess_pr}
  -\delta_{\bs{V}} \nabla_{\bs{B}} \log (\pi(\bs{B}; \beta))
  = \rho \bs{W} \bs{V} \bs{W}
  + \beta (\bs{B} - \hat{\mu}\bs{I}_{d_{\xi}})^{-1} \bs{V} (\bs{B} - \hat{\mu}\bs{I}_{d_{\xi}})^{-1} \\
  + \beta (\hat{L}\bs{I}_{d_{\xi}} - \bs{B})^{-1} \bs{V} (\hat{L}\bs{I}_{d_{\xi}} - \bs{B})^{-1}.
\end{multline}
Finally, substituting \eqref{eq:dir_hess_lkl} and \eqref{eq:dir_hess_pr} into \eqref{eq:dir_hess_post_}, we obtain the directional derivative of the gradient of the negative log posterior,
\begin{multline} \label{eq:dir_hess_post}
  \delta_{\bs{V}} \nabla_{\bs{B}} \cl{L}_k(\bs{B}; \beta)
    = \frac{1}{\nu} \sum_{\ell=1}^k \bs{P}_\ell \bs{V} \bs{s}_\ell \otimes \bs{s}_\ell
    + \rho \bs{W} \bs{V} \bs{W} \\
    + \beta (\bs{B} - \hat{\mu}\bs{I}_{d_{\xi}})^{-1} \bs{V} (\bs{B} - \hat{\mu}\bs{I}_{d_{\xi}})^{-1}
    + \beta (\hat{L}\bs{I}_{d_{\xi}} - \bs{B})^{-1} \bs{V} (\hat{L}\bs{I}_{d_{\xi}} - \bs{B})^{-1}.
\end{multline}

\subsection{Hessian approximation inverse}
To use the Newton's method to solve \eqref{eq:subproblem}, one needs to either solve the linear system $\bs{B}_k \bs{d}_k = \nabla_{\bs{\xi}}F(\bs{\xi}_k)$, or find the inverse of the Hessian approximation $\bs{H}_k = \bs{B}_k^{-1}$, thus $\bs{d}_k = \bs{H}_k \nabla_{\bs{\xi}}F(\bs{\xi}_k)$.

To avoid inverting $\hat{\bs{B}}_k$, we propose to use the Newton--Raphson method to find an approximation $\hat{\bs{H}}_k$ such that $\|\hat{\bs{B}}_k \hat{\bs{H}}_k - \bs{I}_{d_{\xi}} \|_F \le tol_{\bs{H}}$; we can start the Newton--Raphson from the previous Hessian inverse approximation, thus saving a significant computational time.
The Newton--Raphson method for this problem is given by
\begin{equation}
  \hat{\bs{H}}_k^{i+1} = \hat{\bs{H}}_k^i - (R'(\hat{\bs{H}}_k^i))^{-1} R(\hat{\bs{H}}_k^i),
\end{equation}
where $R(\bs{H}) = \bs{B}_k\bs{H} - \bs{I}_{d_{\xi}}$ and $R'(\bs{H}) = \bs{B}_k$, thus,
\begin{equation}
  \hat{\bs{H}}_k^{i+1} = (2 \bs{I}_{d_{\xi}} - \hat{\bs{H}}_k^i \bs{B}_k) \hat{\bs{H}}_k^i.
\end{equation}
The recursive procedure is repeated until $\|R(\hat{\bs{H}}_k^i)\|_F<tol_{\bs{\bs{H}}}$.
Then, we set $\hat{\bs{H}}_k \coloneqq \hat{\bs{H}}_k^i$ as an approximation of $\bs{H}_k$ in optimization.
The convergence of the Newton--Raphson method depends on the initial guess of the inverse, $\hat{\bs{H}}_k^0$; the method converges if the spectral norm of $R(\hat{\bs{H}}_k^0)$ is less than one \cite{pan1985efficient}.
We start from the inverse estimate in the previous iteration, $\hat{\bs{H}}_k^0 = \hat{\bs{H}}_{k-1}$.
If the method does not converge, that is, if the residue increases, we reset $\hat{\bs{H}}_k^i = 2/(\hat{L}+\hat{\mu})\,\bs{I}_{d_{\xi}}$.
This starting point is similar to the one suggested by \cite{pan1985efficient} for Hermitian positive definite $\hat{\bs{B}}_k$, $\hat{\bs{H}}_k^i = 2/(\lambda_{max}+\lambda_{min})\,\bs{I}_{d_{\xi}}$, where $\lambda_{max}$ and $\lambda_{min}$ are, respectively, the largest and the smallest eigenvalues of $\hat{\bs{B}}_k$.


\section{Preconditioned Stochastic Gradient Descent methods}

Instead of the common update of SGD as presented in \eqref{eq:sgd}, we investigate the case where the stochastic gradient estimate is preconditioned with our Hessian inverse approximation,
\begin{equation}\label{eq:precond_sgd}
\xikp = \bs{\xi}_{k} - \eta_k \hat{\bs{H}}_k \bs{\upsilon}_k.
\end{equation}
If the function to be minimized is strongly convex, SGD with fixed step size converges linearly to a neighborhood of the optimum \cite{gower2019sgd}.

Here, we expect the preconditioned SGD to converge to a neighborhood of the optimum faster than vanilla SGD.
One possible side effect of the preconditioning of the stochastic gradient estimate is that the noise in the gradient might be amplified \cite{li2017preconditioned}.
Therefore, it is fundamental to control the error in the gradient estimates using variance control techniques,
e.g., MICE \cite{MICE}, SVRG \cite{johnson2013accelerating}, SARAH \cite{nguyen2017sarah}.
Note that the variance control methods mentioned before are specializations of \eqref{eq:sgd}, only differing in the way that they compute $\bs{\upsilon}_k$.
Moreover, acknowledging the possibly prohibitive cost of updating the Hessian estimate every iteration, we update the Hessian after a given number of iterations or gradient evaluations.

Some variance control methods are based on computing differences between gradients at known points and are thus well-suited to be combined with our approach.
The MICE algorithm \cite{MICE} builds a hierarchy of iterates, then progressively computes gradient differences between subsequent iterates in this hierarchy to control the gradient estimator error.
Since both the differences between the iterates and between the gradients are available, curvature pairs $\bs{s}$ and $\bar{\bs{y}}$ can be obtained without any extra cost in terms of memory or gradient evaluations.
If we refer to the algorithm in \eqref{eq:sgd} using MICE to obtain $\bs{\upsilon}_k$ as SGD-MICE, preconditioning the gradient estimates by our inverse Hessian approximation as in \eqref{eq:precond_sgd} furnishes the SGD-MICE-Bay method.
For the finite sum case, other control variates methods are available, e.g., SVRG\cite{johnson2013accelerating}, SARAH \cite{nguyen2017sarah}, SAGA \cite{defazio2014saga}, SAG \cite{schmidt2017minimizing}.
However, SVRG and SARAH have an advantage over SAGA and SAG here in that curvature pairs can be obtained directly from the algorithms without extra memory or processing costs.
We refer to SVRG and SARAH preconditioned by our inverse Hessian approximation as SVRG-Bay and SARAH-Bay, respectively.

\subsection{Convergence of preconditioned stochastic gradient descent}

Assume that the relative statistical error of the gradient estimator is uniformly bounded above by a constant $\epsilon$ and that an arbitrary lower-bound $\tilde{\mu}$ is imposed on the eigenvalues of the Hessian approximation, i.e., $\tilde{\mu} = \varepsilon \mu$ with $\varepsilon \ge 1$.
On the one hand, a large value of $\varepsilon$ decreases the noise amplification, as previously discussed, and at the same time improves the stability of the preconditioned stochastic gradient descent.
On the other hand, if $\mu$ is indeed the smallest eigenvalue of the true Hessian matrix of $F$, imposing the constraint that the Hessian approximation must have a smallest eigenvalue $\tilde{\mu} \gg \mu$ causes the Hessian approximation to differ from the true Hessian.
In our analysis, we seek to consider the effect of the choice of parameters $\varepsilon$, $\epsilon$, and the step size $\eta$.

\begin{ass}[Regularity and convexity of the objective function] \label{ass:regularity}
The objective function $F$ has Lipschitz continuous second-order derivatives and is $\mu$-convex. Namely, for all $\bs{x}, \bs{y} \in \Xi$,
\begin{align}
  \| \hess{\bs{x}} - \hess{\bs{y}} \| &\le M \norm{\bs{x} - \bs{y}} \\
  F(\bs{y}) &\ge F(\bs{x}) + \grad{\bs{x}} \cdot (\bs{y} - \bs{x}) + \frac{\mu}{2} \norm{\bs{y} - \bs{x}}^2.
\end{align}
\end{ass}

\begin{ass}[Unbiasedness and bounded statistical error of gradient estimator]
\label{ass:gradient}
The gradient estimator $\bs{\upsilon}_k$ is an unbiased estimator of $\grad{\xik}$ with an upper bound $\epsilon$ on its relative statistical error,
\begin{align}
  \Expc{\bs{\upsilon}_k}{\xik} &= \grad{\xik} \\
  \label{eq:grad_condition}
  \Expc{\norm{\bs{\upsilon}_k - \Expc{\bs{\upsilon}_k}{\xik}}^2}{\xik} &\le  \epsilon^2 \norm{\grad{\xik}}^2.
\end{align}
\end{ass}

\begin{ass}[Bounded error of Hessian approximation]
  \label{ass:hessian}
  The approximation $\Bk$ of $\hess{\xik}$ has a lower-bound $\tilde{\mu} \coloneqq \varepsilon \mu$ on its eigenvalues with $\varepsilon \ge 1$ and the product of the true Hessian by the approximation inverse $\Hk$ results in the identity matrix added by an error that is bounded above in norm by a linear function of $\varepsilon$,
  \begin{align}
    \Bk &\succeq \tilde{\mu}\bs{I}_{\dimxi} \coloneqq \varepsilon \mu, \, \varepsilon \ge 1, \\
    \hess{\xik} \Hk & = \bs{I} + \cl{E}_k, \\
    \norm{\cl{E}_k} &\le C_{\varepsilon} (\varepsilon - 1),
  \end{align}
  where $C_{\varepsilon}>0$ is a constant.
\end{ass}

\begin{thm}[Convergence of preconditioned stochastic gradient descent]
  \label{thm:conv_sgd}
  If Assumptions \ref{ass:regularity}, \ref{ass:gradient}, and \ref{ass:hessian} are satisfied and
  \begin{equation} \label{eq:cond_conv}
      \norm{\grad{\xik}}
     <
     \frac{2 \varepsilon^2 \mu^2 \big(
      1 - \epsilon - C_{\varepsilon} (1 + \epsilon) (\varepsilon - 1)
      - \frac{r}{\eta}
     \big)}{
     M \eta (1 + \epsilon^2)
     }
  \end{equation}
  holds for all $k>0$,
  then the gradients of the iterates generated by the preconditioned stochastic gradient descent as described in \eqref{eq:precond_sgd} converge linearly wirh rate $1-r$,
  \begin{equation}
    \norm{\grad{\xikp}}
    \le
    (1 - r)^{k+1} \norm{\grad{\bs{\xi}_0}}.
  \end{equation}
\end{thm}

\begin{proof}
From Assumption \ref{ass:regularity}, we obtain an upper bound for the norm of the gradient at $\xikp$,
\begin{align}
  \norm{\grad{\xikp}}
  &\le \norm{\grad{\xik} + \hess{\xik} (\xikp - \xik)} + \frac{M}{2} \norm{\xikp - \xik}^2.
\end{align}
Using the preconditioned SGD update in \eqref{eq:precond_sgd}, one can further develop the equation above as
\begin{align}
  \norm{\grad{\xikp}}
  & \le \norm{\grad{\xik} - \eta \hess{\xik} (\Hk \bs{\upsilon}_k)}
   + \frac{M \eta^2}{2} \norm{\Hk \bs{\upsilon}_k}^2 \nonumber \\
  & \le 
  \!\begin{multlined}[t][10cm]
  \norm{\grad{\xik} - \eta \hess{\xik} \Hk \grad{\xik }} \\
    + \norm{\eta \hess{\xik} \Hk (\bs{\upsilon}_k - \grad{\xik})}
    + \frac{M \eta^2}{2} \norm{\Hk \bs{\upsilon}_k}^2
  \end{multlined}
  \\
  & \le 
  \!\begin{multlined}[t][10cm]
  \norm{\bs{I}_{\dimxi} - \eta \hess{\xik} \Hk}
      \norm{\grad{\xik}} \\
    + \eta \norm{\hess{\xik} \Hk} \norm{\bs{\upsilon}_k - \grad{\xik}}
    + \frac{M \eta^2}{2} \norm{\Hk}^2 \norm{\bs{\upsilon}_k}^2.
  \end{multlined}
\end{align}
Taking the expectation conditioned on $\xik$ and using Assumption \ref{ass:gradient},
\begin{multline}
  \Expc{\norm{\grad{\xikp}}}{\xik}
  \le
  \norm{\bs{I}_{\dimxi} - \eta \hess{\xik} \Hk}
      \norm{\grad{\xik}} \\
  + \eta \epsilon \norm{\hess{\xik} \Hk} \norm{\grad{\xik}} 
  + \frac{M \eta^2}{2} (1 + \epsilon^2) \norm{\Hk}^2 \norm{\grad{\xik}}^2.
\end{multline}
Using Assumption \ref{ass:hessian},
\begin{align}
  \Expc{\norm{\grad{\xikp}}}{\xik}
  & \le
  \!\begin{multlined}[t][10cm]
  \left(
  \norm{\bs{I}_{\dimxi} - \eta (\bs{I}_{\dimxi} + \cl{E}_k)}
  + \eta \epsilon \norm{\bs{I}_{\dimxi} + \cl{E}_k} \right)
  \norm{\grad{\xik}} \\
  + \frac{M \eta^2}{2 \tilde{\mu}^2} (1 + \epsilon^2) \norm{\grad{\xik}}^2
  \end{multlined}
  \nonumber \\
  & \le
  \!\begin{multlined}[t][10cm]
  \big(
  (1 - \eta) +
  \eta \norm{\cl{E}_k}
  + \eta \epsilon (1 + \norm{\cl{E}_k})
  \big)
  \norm{\grad{\xik}} \\
  + \frac{M \eta^2}{2 \tilde{\mu}^2} (1 + \epsilon^2) \norm{\grad{\xik}}^2 
  \end{multlined}
  \nonumber \\
  & \le
  \!\begin{multlined}[t][10cm]
  \big(
  1 - \eta (1 - \epsilon) + \eta (1 + \epsilon) C_{\varepsilon} (\varepsilon - 1)
  \big) \norm{\grad{\xik}} \\
  + \frac{M \eta^2}{2 \varepsilon^2 \mu^2} (1 + \epsilon^2) \norm{\grad{\xik}}^2.
  \end{multlined}
  \label{eq:sgd_proof_last}
\end{align}
Letting the right hand side of the above inequality be strictly less than $(1 - r) \norm{\grad{\xik}}$ results in \eqref{eq:cond_conv}.
Taking expectation and unrolling the recursion complete the proof.

\end{proof}
Theorem \ref{thm:conv_sgd} shows that increasing the parameter $\varepsilon$ that controls the smallest eigenvalue of the Hessian approximation improves the stability of preconditioned stochastic gradient descent.
Moreover, controlling $\varepsilon$ and the step size $\eta$ to satisfy \eqref{eq:cond_conv} for a given $r$ guarantees global convergence.
By letting $\epsilon=0$, thus resulting in deterministic gradient observations; $\varepsilon=1$, meaning that $\tilde{\mu}=\mu$; and unitary step size $\eta=1$, one recovers the quadratic convergence of the Newton method with the condition that $\norm{\grad{\xik}} < 2\mu^2/M$ (see \eqref{eq:sgd_proof_last}).

\begin{figure}[h]
  \centering
  \includegraphics[width=.5\linewidth]{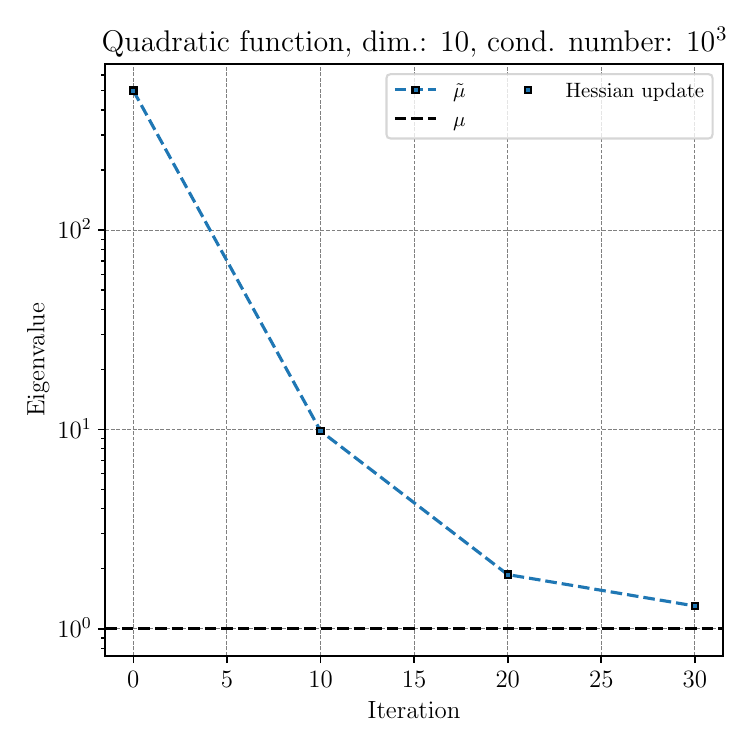}
  \caption{Value of $\tilde{\mu}$ versus iteration for SGD-MICE preconditioned with our Hessian inverse approximation applied to the minimization of the quadratic problem in example in Section \ref{sec:quadratic}.}
  \label{fig:min_eig}
\end{figure}

Since our analysis depends on quantities that are hard to characterize in practice, we do not use it to control $\varepsilon$ and $\eta$.
Instead, we start optimization with $\varepsilon=(L + \mu)/(2\mu)$ and control the parameter $\rho$ of the Frobenius regularization of the prior; the regularization term penalizes large changes in the Hessian approximation. Thus, at every Hessian approximation update, $\varepsilon$ progressively decays.
Figure \ref{fig:min_eig} illustrates this behavior for the Example in Section \ref{sec:quadratic}.

\subsection{Time complexity}

Given that the interval between Hessian updates is set to be $\bigO{d_{\bs{\xi}}}$, then the extra cost of finding the Hessian approximation using our Bayesian approach is $\bigO{d_{\bs{\xi}}^2}$.
The number of Newton iterations and conjugate gradient iterations do not depend on the problem's dimensionality but on the tolerance set for each algorithm.
The cost of each conjugate gradient iteration is the cost of evaluating $\delta_{\bs{V}} \nabla_{\bs{B}} \cl{L}_k(\bs{B}, \beta)$ in \eqref{eq:dir_hess_post}, which is $\bigO{d_{\bs{\xi}}^3}$ due to matrix multiplications.
If, instead of using the conjugate gradient method to compute the Newton direction, one builds the fourth-order tensor of second-order derivatives of the log posterior with respect to the Hessian, the time complexity increases to $\bigO{d_{\bs{\xi}}^4}$.

In each iteration, we have to compute $P_k$, the precision matrix for each curvature pair.
If a full matrix is computed, we have a complexity of $\bigO{d_{\bs{\xi}}^2}$ each iteration.
If, however, we opt to ignore the off-diagonal terms, this complexity is improved to $\bigO{d_{\bs{\xi}}}$.
We can further decrease the memory overhead of our method by storing only a scalar, the trace of the covariance matrix, for each curvature pair.
The reasoning behind this approach is that the trace of the covariance matrix works as a weight for the curvature pairs; in the curvature pairs where the gradient difference is better estimated, the respective curvature pair has a larger weight when computing the likelihood.
When our Bayesian Hessian approximation is used with the MICE gradient estimator, the trace of the covariance matrix is available without any extra cost.
In numerical tests, we observed that keeping the trace of the covariance matrix resulted in the best option regarding the balance between accuracy and computational time.

We do not discuss here the cost to compute $\bs{\upsilon}_k$ satisfying condition \eqref{eq:grad_condition}, as it depends on the method used to estimate it.
For example, if $\bs{\upsilon}_k$ is a Monte Carlo estimator with adaptive batch sizes, the cost of evaluating $\bs{\upsilon}_k$ is $\bigO{\epsilon^{-2} \norm{\grad{\xik}}^{-2}}$.
If the conditions of Theorem 1 are satisfied, the gradient sampling cost is $\bigO{\epsilon^{-2} (1-r)^{-2k}}$.

\section{Numerical examples}
\label{sec:numerics}

In this section, we solve some numerical examples with and without our Bayesian approach to validate its performance.
In all cases, we present the convergence of the optimality gap, $F(\bs{\xi}_k) - F(\bs{\xi}^*)$, and in some cases, the convergence of the central-path Newton-CG used to find the Hessian approximations.
The parameters of the Bayesian approach are kept fixed for all the numerical examples, namely,
the logarithmic barrier parameter, $\beta=10^{-2}$; the Frobenius norm regularizer parameter, $\rho=10^{-2}$; the tolerance on the norm of the gradient of the negative log posterior, $tol=10^{-6}$; the constraint on the lower bound of the eigenvalues of $\bs{B}$, $\tilde{\mu}$, is set as the strong-convexity parameter $\mu$; the relaxation parameter,
$\alpha=1.05$; the number of central-path steps is six; and the factor of decrease of both $tol_i$ and $\beta_i$ for each $i$-th central-path step is $\gamma=2$.
Moreover, to keep our Bayesian approach competitive in terms of computational time and memory allocation, we do not compute the full matrix $P_\ell$ for each curvature pair, using the inverse of the trace of the covariance matrix instead.
Note that the trace of the covariance matrix of the difference between gradients is available from MICE since it is used to control the relative error of the gradient \cite{MICE}.

The results presented in this section are obtained using a python implementation of the proposed method whose source code can be found at GitHub\footnote{\url{https://github.com/agcarlon/bayhess}}.
Moreover, the source code of the MICE gradient estimator, also used in this section, can be found at Bitbucket\footnote{\url{https://bitbucket.org/agcarlon/mice}}.
A repository with all the code necessary to reproduce the results presented in this section is available at GitHub\footnote{\url{https://github.com/agcarlon/bayhess_numerics}}.

\subsection{Quadratic function} \label{sec:quadratic}

The first numerical test is a quadratic function with noise added to its gradient.
We want to minimize $\bb{E}_{\bs{\theta}} [f(\bs{\xi}, \bs{\theta}) | \bs{\xi}]$ with
\begin{equation}
 f(\bs{\xi}, \bs{\theta}) = \frac{1}{2} \bs{\xi} \cdot \bs{A} \, \bs{\xi}
 - \bs{b} \cdot \bs{\xi} + \bs{\xi} \cdot \sqrt{\bs{\Sigma}} \, \bs{\theta},
\end{equation}
where $\bs{A}$ is a square symmetric positive-definite matrix with $1 \preceq \bs{A} \preceq \kappa$, $\bs{\Sigma}$ is a symmetric and positive definite covariance matrix, $\bs{\theta}$ is a vector of independent standard normal distributed random variables, and the vector $\bs{b}$ is defined as $\bs{b} = \bs{1}_{d_{\bs{\xi}}}$, where $\bs{1}_{d_{\bs{\xi}}}$ is a vector of ones of dimension $d_{\bs{\xi}}$.
The eigenvalues of $\bs{A}$ are sampled uniformly between $1$ and $\kappa$, with necessarily one eigenvalue being $1$ and another being $\kappa$.
The gradient of $f$ with respect to $\bs{\xi}$ is
\begin{equation}
  \nabla f(\bs{\xi}, \bs{\theta}) = \bs{A} \, \bs{\xi} - \bs{b} + \sqrt{\bs{\Sigma}} \, \bs{\theta},
\end{equation}
thus, the covariance matrix of $\nabla f$ is $\bs{\Sigma}$.

We solve this problem for two different values of the condition number $\kappa$, namely $10^3$ and $10^6$, and for a fixed number of dimensions $d_{\bs{\xi}}=10$.
We compare the performances of different methods in the solution of this problem: SGD, SGD-MICE, SGD-Bay with decreasing step size, SGD-Bay with fixed step size, and SGD-MICE-Bay.
In all cases, we run the optimization until $10^7$ gradients are evaluated.
For SGD-Bay, since the number of gradient evaluations per iteration is fixed, we know the number of iterations a priori and thus choose the interval between Hessian updates to have $15$ equally distributed updates.
For SGD-MICE-Bay, we opt to have an update every $d_{\bs{\xi}}$ iterations.
The step sizes and batch sizes used for each method are presented in Table \ref{tab:quadratic}, where, for the cases when MICE is used, we show the tolerance $\epsilon$ on the relative statistical error of the gradient.

\begin{table}
  \centering
 \renewcommand{\arraystretch}{1.5}
 \caption{Step sizes and batch sizes used for Example \ref{sec:quadratic}. Here, $L$ is the Lipschitz constant of $\nabla f$, $\mu$ is the strong-convexity constant, $k$ is the iteration, and $\epsilon$ is the tolerance on the relative statistical error.}
 \label{tab:quadratic}
 \begin{tabular}{lcc}
  \toprule
  Method          & Step size                           & Batch size     \\ \midrule
  SGD             & $\frac{1}{L \sqrt{k}}$              & $1$            \\
  SGD-MICE        & $\frac{2}{(L + \mu)(1+\epsilon^2)}$ & $\epsilon=1$   \\
  SGD-Bay (decr.) & $\frac{1}{\sqrt{k}}$                & $10$           \\
  SGD-Bay (fixed) & $1$                                 & $1000$         \\
  SGD-MICE-Bay    & $1/(1+\epsilon^2)$                  & $\epsilon=0.5$ \\
  \bottomrule
 \end{tabular}
\end{table}

\begin{figure}[h]
  \centering
 \includegraphics[width=.48\linewidth]{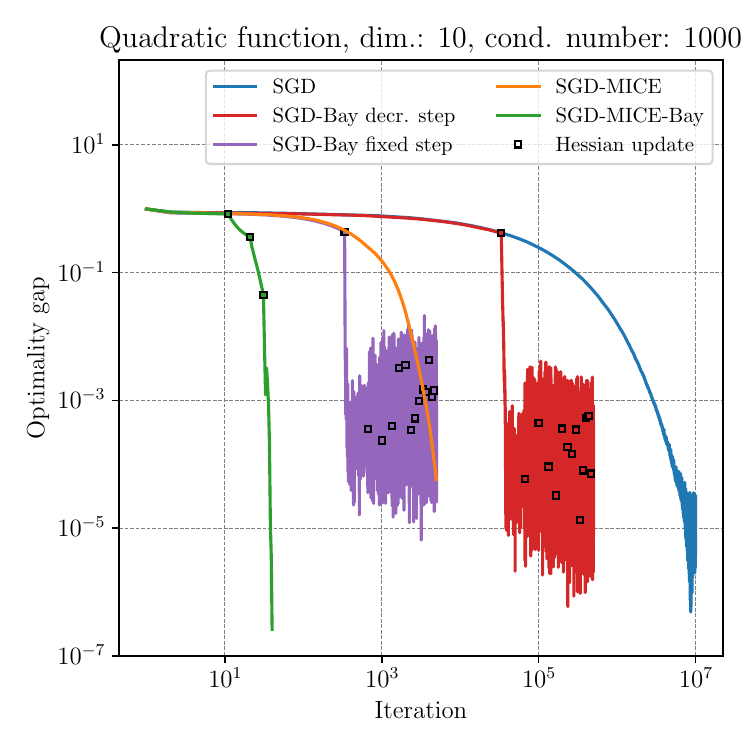}
 \includegraphics[width=.48\linewidth]{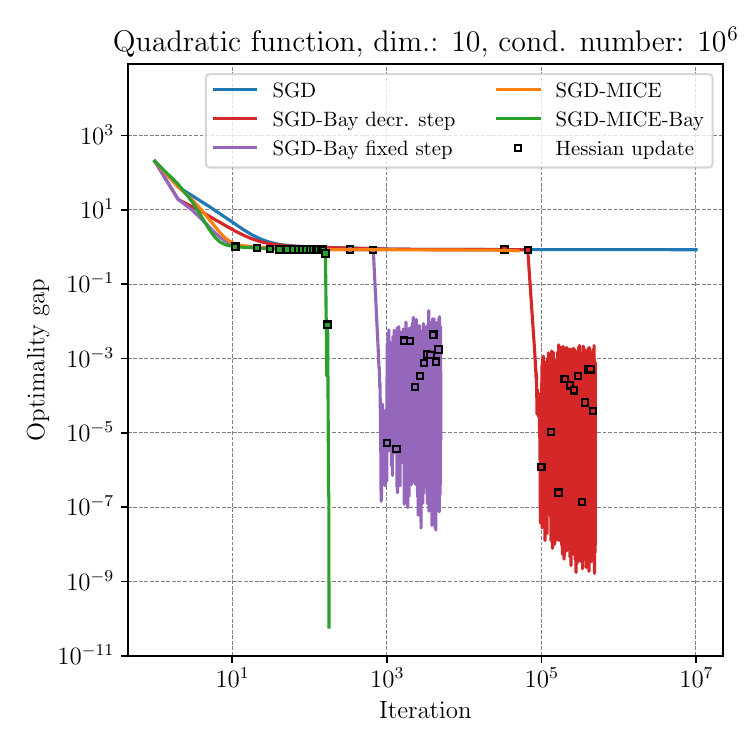}
 \caption{Convergence of the optimality gap versus iteration with $\kappa=1000$ (left) and $\kappa=10^6$ (right) for Example \ref{sec:quadratic}.
 The Hessian updates are marked with colored squares with black edges.
 In both cases, SGD-MICE-Bay reduced the optimality gap more than the other methods and in a much smaller number of iterations.
 Moreover, SGD-MICE-Bay did not oscillate around the optimum like both SGD-Bay cases.}
 \label{fig:ex1_conv}
\end{figure}

In Figure \ref{fig:ex1_conv}, we present the optimality gap versus iteration for both $\kappa=10^3$ (left) and $\kappa=10^6$ (right) for the previously mentioned methods.
In the cases where our Bayesian approach is used, the Hessian updates are illustrated as colored squares with black edges.
For the case where $\kappa=1000$, SGD-MICE-Bay achieved, in $39$ iterations, a lower optimality gap than SGD in $10^7$ iterations.
Moreover, SGD presents oscillations around the optimum, making it difficult to have a reliable optimum estimate.
In the case of SGD-Bay, with and without the decreasing step size, we observe oscillations of larger amplitude due to the preconditioning with our Hessian approach.
This behavior is due to the lack of variance control of SGD; not even the decreasing step size is enough to avoid the oscillations.
In the case of an even larger condition number, $\kappa=10^6$, SGD-MICE-Bay far outperforms the other methods while requiring $178$ iterations.
Both SGD and SGD-MICE, being first-order methods, cannot improve the optimality gap significantly once close to the optimum.
As for SGD-Bay, the same oscillating behavior is observed for both the fixed and decreasing step cases.
From these results, it is clear that, in large condition numbers cases, using second-order information is necessary for an efficient convergence.
Moreover, we observe that controlling the variance of the gradient estimator is of central importance to maintaining convergence in the stochastic setting, avoiding the amplification of noise due to the preconditioning.

One crucial aspect of our method is the cost to update the Hessian approximation, i.e., the cost of the Newton-CG method in this setting.
For SGD-MICE-Bay in the $\kappa=1000$ case, we present the Newton-CG convergence of the gradient norm $\|\nabla_{\bs{B}} \cl{L}\|$ per Newton iteration.
The black numbers on top of each update denote the number of conjugate gradient steps needed to find the Newton direction.
Each color represents a different central-path step, i.e., at each color change, both $\beta$ and $tol$ are decreased by a factor $\gamma=2$, with each $tol$ denoted by a dash-dotted line.
For this problem, the runtime of each of the Hessian updates is at fractions of seconds.
Given that $d_{\bs{\xi}}=10$, the Hessian matrix has $100$ components.
Still, in the worst case, $51$ conjugate gradient iterations were needed to find the Newton direction, showing that the problem of finding $\hat{\bs{B}}$ is well-conditioned.
Moreover, the total number of Newton iterations in the three cases did not exceed $25$.

\begin{figure}[h]
  \centering
\includegraphics[width=.7\linewidth]{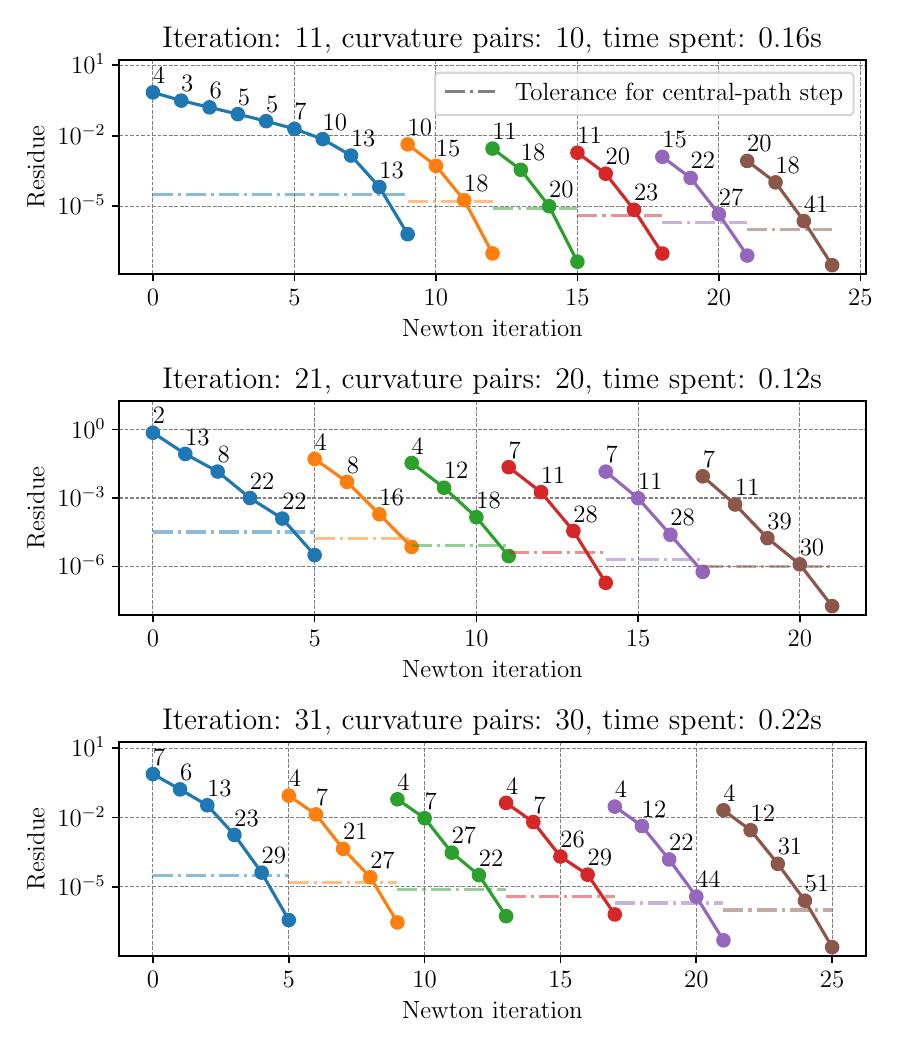}
\caption{For each of the three Hessian updates of SGD-MICE-Bay with $\kappa=1000$ in Example \ref{sec:quadratic}, we present the convergence of the gradient norm for the Newton-CG method.
The number of curvature pairs available and the time taken are presented for each Hessian update.
Each color represents a different step of the central-path method with a different logarithmic barrier parameter $\beta$ and a different residue tolerance.
The tolerance for each central-path step is presented as a dash-dotted line.
The number above each Newton iteration represents the number of conjugate gradient iterations needed to find the Newton direction.
}
\label{fig:ex1_newton_cg}
\end{figure}

\begin{figure}
 \includegraphics[width=.48\linewidth]{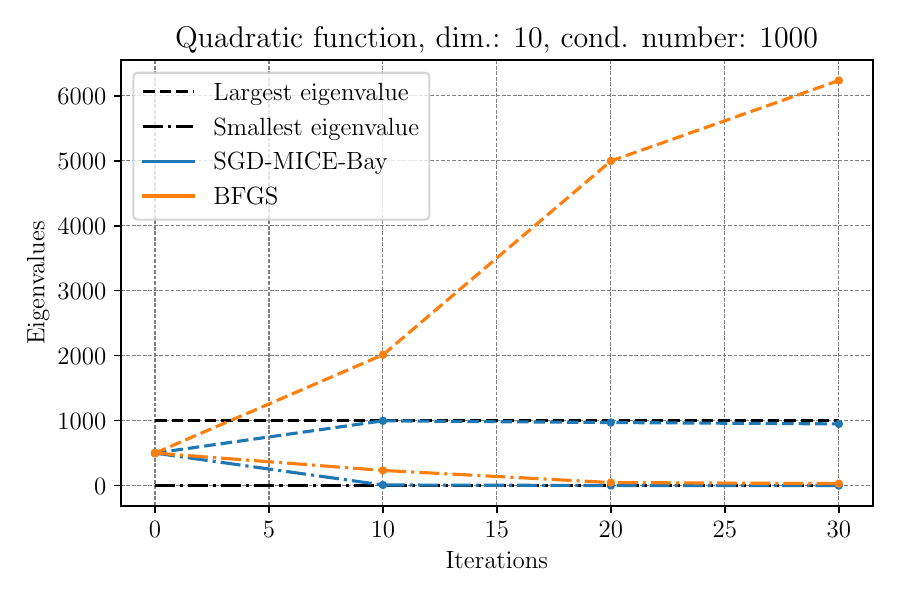}
 \includegraphics[width=.48\linewidth]{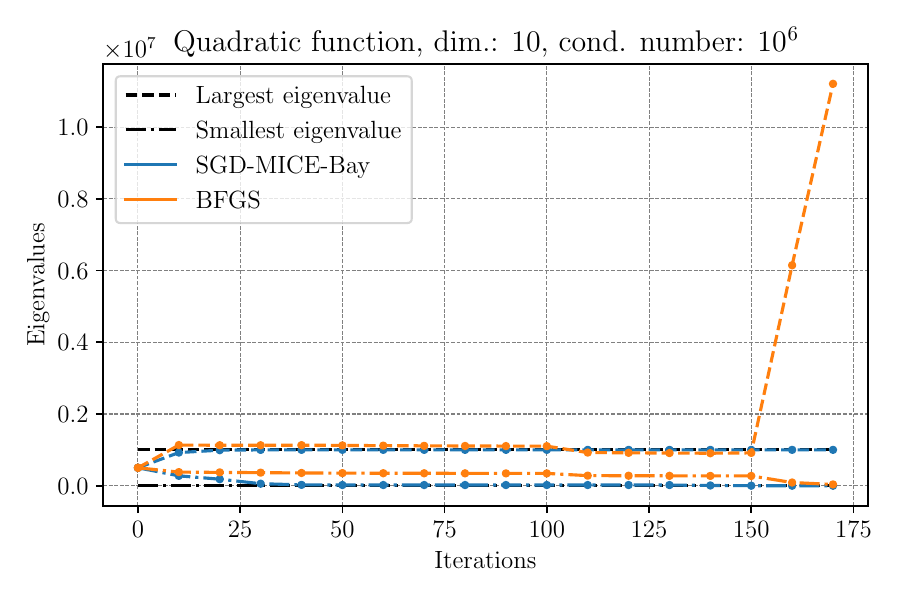}
 \caption{Maximum and minimum eigenvalues for both our Bayesian approach and BFGS for Example \ref{sec:quadratic} with $\kappa=1000$ (left) and $\kappa=10^6$ (right).
 We use the same curvature pairs for both cases, obtained from SGD-MICE-Bay.
 The true extremes of the eigenvalues are presented as black lines.
 Dashed lines represent the largest eigenvalues, and dash-dotted lines represent the smallest eigenvalues.
 It is clear that our approach gets closer to the actual extreme eigenvalues than BFGS, which is farther from the smallest eigenvalue than our method and exceeds the value of the largest eigenvalue.
 }
 \label{fig:ex1_err_and_eigs}
\end{figure}

Controlling the eigenvalues maxima and minima is essential to keep the stability and efficiency of preconditioned SGD.
To motivate why our approach is better suited to stochastic optimization than using the BFGS formulae to approximate the Hessian, we, for every Hessian update of SGD-MICE-Bay, also compute a Hessian approximation using BFGS.
In Figure \ref{fig:ex1_err_and_eigs}, we present, for both $\kappa=1000$ (left) and $\kappa=10^6$ (right), the extreme eigenvalues for our Bayesian approach, for BFGS, and the true extremes.
The same curvature pairs are used for both our Bayesian approach and BFGS.
For both condition number cases, the BFGS method can keep the minimum eigenvalue above the actual value; however, the largest eigenvalue far exceeds the true largest eigenvalue.
In contrast, our Bayesian approach keeps the eigenvalues of the Hessian approximation between the extreme values due to the logarithmic barrier constraints.

\begin{figure}
  \centering
 \includegraphics[width=.5\linewidth]{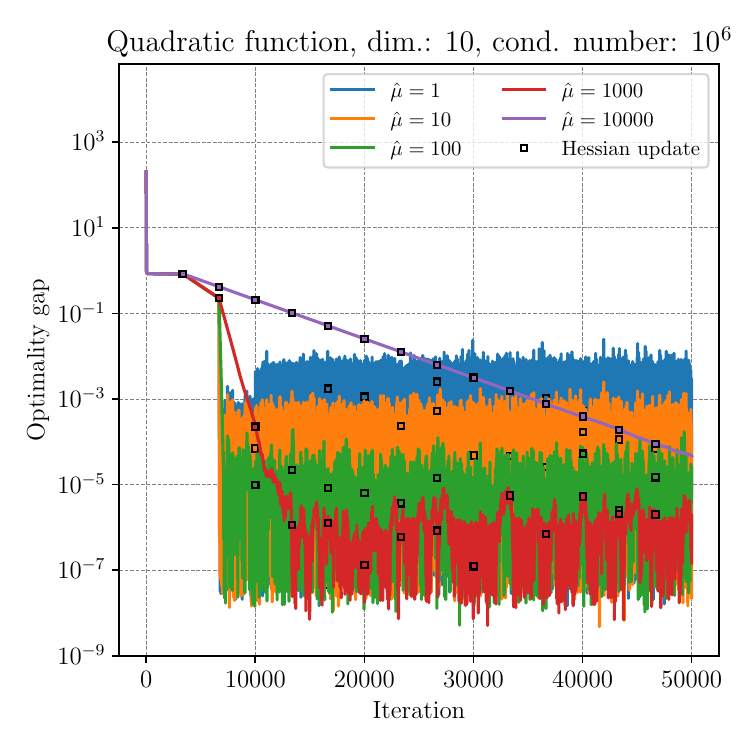}
 \caption{
 Optimality gap versus iteration for SGD-Bay with fixed step size on Example \ref{sec:quadratic} for different values of $\hat{\mu}$.
 The colored squares with black edges represent the Hessian updates of each method.
 The different values of $\hat{\mu}$ result in different amplitudes for the oscillations around the optimum of SGD-Bay; larger values of $\hat{\mu}$ result in smaller amplitudes of oscillations (cf. \eqref{eq:noise}).
 }
 \label{fig:ex1_stab_test}
\end{figure}

As motivated in \eqref{eq:noise} and Theorem \ref{thm:conv_sgd}, increasing $\hat{\mu}$ decreases the noise amplification of preconditioning the gradients by our Bayesian Hessian approximation.
Figure \ref{fig:ex1_stab_test} presents a comparison between SGD-Bay with a fixed step size for different values of $\hat{\mu}$.
One thing to notice is the reduction of the amplitude of the oscillations as $\hat{\mu}$ increases.
Moreover, we observe a clear linear convergence for the more conservative cases of $\hat{\mu}=1000$ and $\hat{\mu}=10000$.
Comparing these results with the ones presented in Figure \ref{fig:ex1_conv} for the same condition number $\kappa=10^6$, it is clear that increasing $\hat{\mu}$ improves the performance of SGD-Bay even without controlling the gradient statistical error or decreasing the step size.
However, in the general convex case, SGD-Bay with a fixed step size will converge to a neighborhood of the optimum and oscillate around it with a fixed amplitude; increasing $\hat{\mu}$ or decreasing the step size may reduce the amplitude of the oscillations but not circumvent them.

\subsection{Logistic regression} \label{sec:logreg}

Training logistic regression models for classification of finite data is an empirical risk minimization problem, thus a problem of minimizing a sum of functions.
When training a model on a large dataset, it is common to use SGD to reduce the iteration cost by subsampling the data points used in each iteration.
Since inexact information about the gradient is used in each iteration, deterministic quasi-Newton methods usually fail to recover the Hessian and its inverse.
Here, we use our Bayesian framework to recover the Hessian from noisy gradient observations and precondition the gradient estimates with the Hessian inverse.
We compare plain SGD, SGD-MICE, SGD-MICE-Bay, SVRG, SVRG-Bay, SARAH, and SARAH-Bay on the task of training a binary classification logistic regression model with $\ell_2$ regularization on five different datasets, \emph{mushrooms}, \emph{ijcnn1}, \emph{w8a} \cite{kohavi1996scaling}, \emph{cod-rna} \cite{uzilov2006detection}, and \emph{HIGGS} \cite{baldi2014searching}\footnote{
 Obtained from LibSVM datasets at
 \url{https://www.csie.ntu.edu.tw/~cjlin/libsvmtools/datasets/binary.html}, except for the HIGGS dataset that was obtained from \url{https://archive.ics.uci.edu} \cite{Dua2022}.
}.

The $\ell_2$-regularized log loss function we seek to minimize is
\begin{equation}\label{eq:logistic}
 F(\bs{\xi}) = \frac{1}{N} \sum_{i=1}^N f\left(\bs{\xi}, \bs{\theta}_i= (\bs{x}_i, y_i)\right) = \frac{1}{N} \sum_{i=1}^N \log(1 + \exp(-y_i \, \bs{\xi} \cdot \bs{x}_i))
 + \frac{\lambda}{2} \|\bs{\xi}\|^2,
\end{equation}
where each data point $(\bs{x}_i, y_i)$ is such that $\bs{x}_i \in \bb{R}^{d_{\bs{\xi}}}$ and $y_i \in \{-1,1\}$.
For all datasets, we use a regularization parameter $\lambda = 10^{-5}$.

For SGD, only a singleton sample is used, i.e., the minibatch size is set to 1.
Also, as in the example in Section \ref{sec:quadratic}, a decreasing schedule is used for the step size.
For SGD-MICE, the sample size is adaptively controlled to keep the relative error of the gradient estimates below $\epsilon=0.5$ with a step size of $\eta=2/((L+\mu)(1+\epsilon^2))$ \cite{MICE}.
For the \emph{HIGGS} dataset, the choice of $\epsilon=0.5$ was too conservative, thus we used $\epsilon=0.8$.
As for SGD-MICE-Bay, we used the same values of $\epsilon$, however, since we precondition the gradient estimates with the inverse of the Hessian, we use a step size of $\eta = (1+\epsilon^2)^{-1}$.
Therefore, with the initial Hessian guess with diagonal entries $B_{ii} = \frac{L+\mu}{2}$, SGD--MICE-Bay recovers exactly SGD-MICE with the optimal step size.
Similarly, for SVRG and SARAH, we use a fixed step size of $\eta_k=0.1/L$, whereas for SVRG-Bay and SARAH-Bay we use $\eta_k=0.1$.
For SVRG, SVRG-Bay, SARAH, and SARAH-Bay, we start with a full batch, iterate using a mini-batch size of five, and restart every two epochs.

In Table \ref{tab:logreg_datasets}, we present the sizes, number of features, and condition numbers of the resulting optimization problems for the five datasets.
The used datasets differ significantly in their characteristics.

\begin{table}[h]
  \centering
 \caption{For each of the datasets used in Example \ref{sec:logreg}, we present its data size, number of features, and condition number of the optimization problem.}
 \label{tab:logreg_datasets}
 \begin{tabular}{l r r r}
  \toprule
  dataset          & \multicolumn{1}{c}{size}            & features         & \multicolumn{1}{c}{$\kappa$}  \\ \midrule
  \emph{mushrooms} & $8124$          & $112$            & $5.25 \times 10^{5}$      \\
  \emph{ijcnn1}    & $49990$         & $22$             & $3.76 \times 10^{4}$      \\
  \emph{w8a}       & $49749$         & $\bs{300}$       & $2.91 \times 10^5$         \\
  \emph{cod-rna}   & $59535$         & $8$              & $\bs{6.48 \times 10^{9}}$ \\
  \emph{HIGGS}     & $\bs{11000000}$ & $28$             & $2.50 \times 10^4$         \\
  \bottomrule
 \end{tabular}
\end{table}

\begin{figure}[h]
  \includegraphics[width=.48\linewidth]{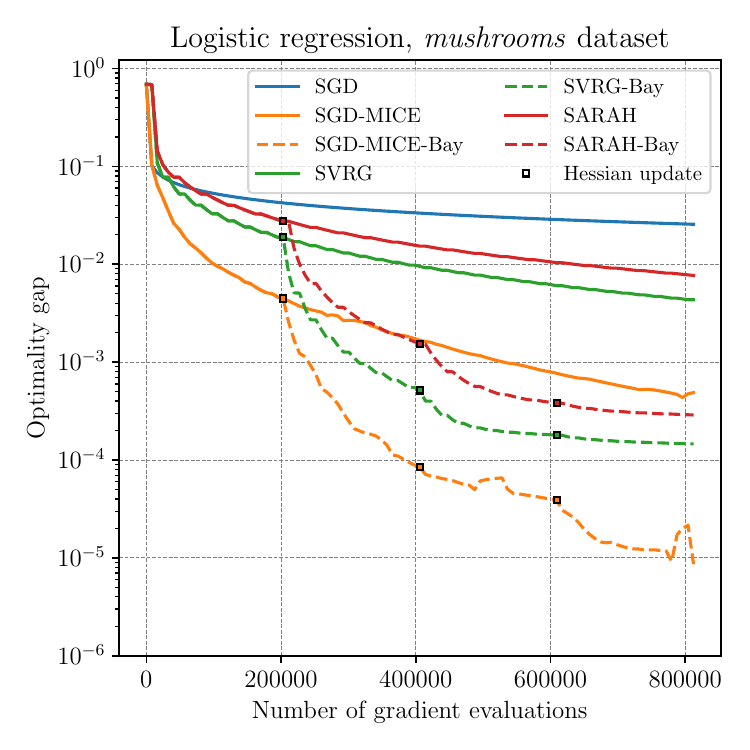}
  \includegraphics[width=.48\linewidth]{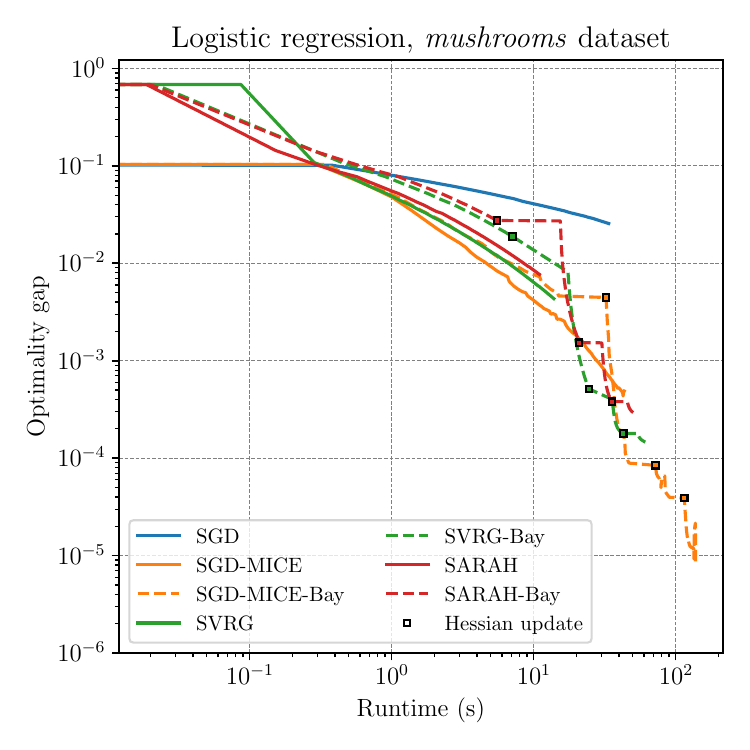}
 \caption{Convergence of the optimality gap versus iteration (left) and runtime (right) for Example \ref{sec:logreg}, \emph{mushrooms} dataset. The Hessian updates are marked with colored squares with black edges.
 }
 \label{fig:ex2_mush_conv}
\end{figure}

In Figure \ref{fig:ex2_mush_conv}, we present the convergence for the optimality gap versus iterations (left) and versus runtime in seconds (right) for the \emph{mushrooms} dataset.
For SGD, SGD-MICE, SVRG, and SARAH, we present their convergences as solid lines and their respective preconditioned cases as dashed lines.
For all the methods, preconditioning the gradient estimates using our Bayesian approach improved the final optimality gap in at least one order of magnitude.
In terms of convergence of the optimality gap versus runtime, we observe that computing the Hessian approximation adds an overhead to the optimization process; however, the improvement in convergence compensates for this overhead.
In this example, SGD-MICE-Bay achieved an optimality gap ten times smaller than that found by SGD-MICE and more than two orders of magnitude smaller than the one found by SGD.
Preconditioning the gradient estimates of SVRG and SARAH with our Bayesian Hessian inverse approximation reduced the final optimality gap found by more than one order of magnitude.

\begin{figure}[h]
  \centering
 \includegraphics[width=.6\linewidth]{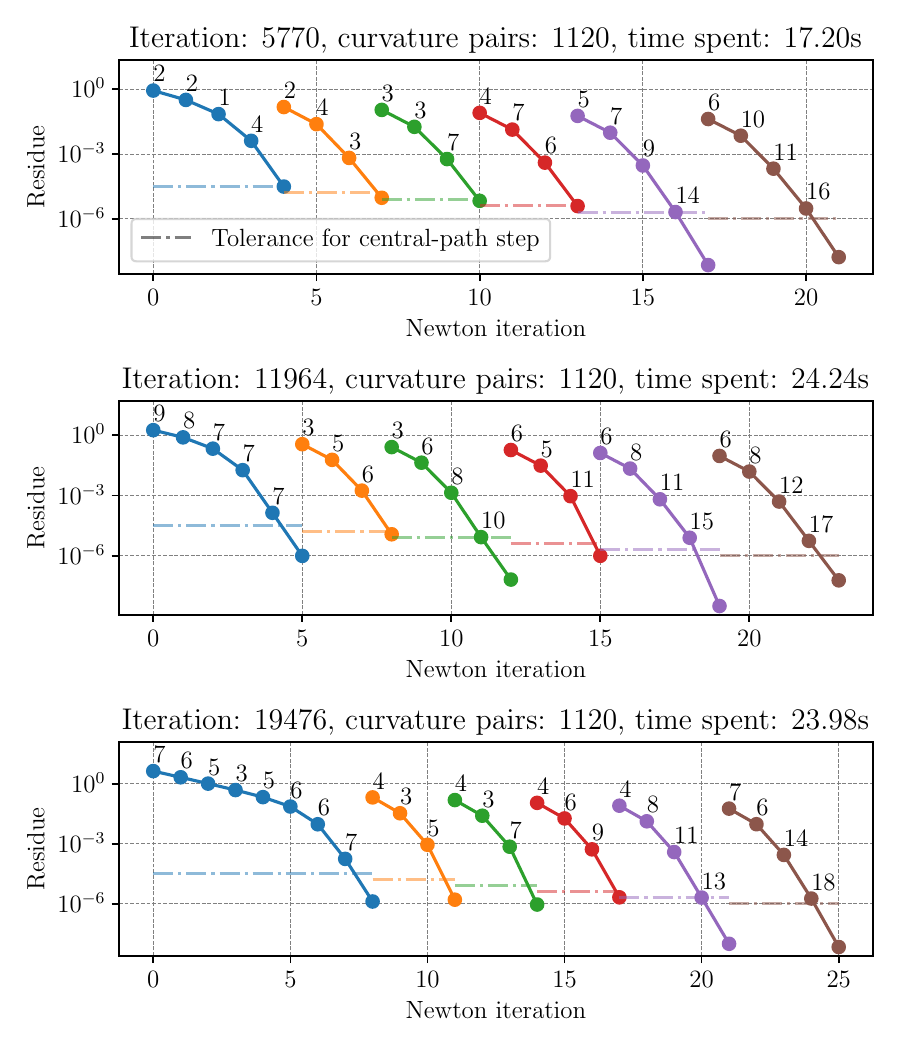}
 \caption{Convergence of the Newton-CG method for each of the Hessian updates of SGD-MICE-Bay in Example \ref{sec:logreg}, \emph{mushrooms} dataset.
 The number of curvature pairs available and the time taken are presented for each Hessian update.
 Each color represents a different step of the central-path method with a different logarithmic barrier parameter $\beta$ and a different residue tolerance.
 The tolerance for each central-path step is presented as a dash-dotted line.
 The number above each Newton iteration represents the number of conjugate gradient iterations needed to find the Newton direction.
 }
 \label{fig:ex2_mush_newton_cg}
\end{figure}

In Figure \ref{fig:ex2_mush_newton_cg}, the convergence of the Newton-CG method is presented for each of the Hessian updates of SGD-MICE-Bay for the \emph{mushrooms} dataset.
Note that a relatively small number of Newton steps are needed for each central-path step; $8$ in the worst case.
Moreover, the number of conjugate gradient steps does not exceed $18$, which is a remarkably small number of iterations considering that we have a problem with $112 \times 112$ variables.
This result indicates that the problem is well-conditioned.

\begin{figure}
  \includegraphics[width=.48\linewidth]{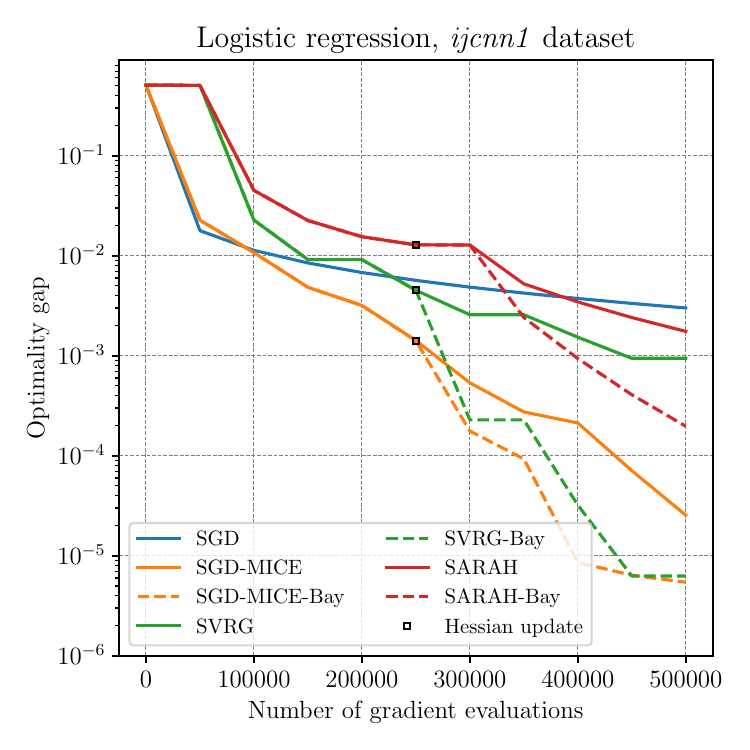}
  \includegraphics[width=.48\linewidth]{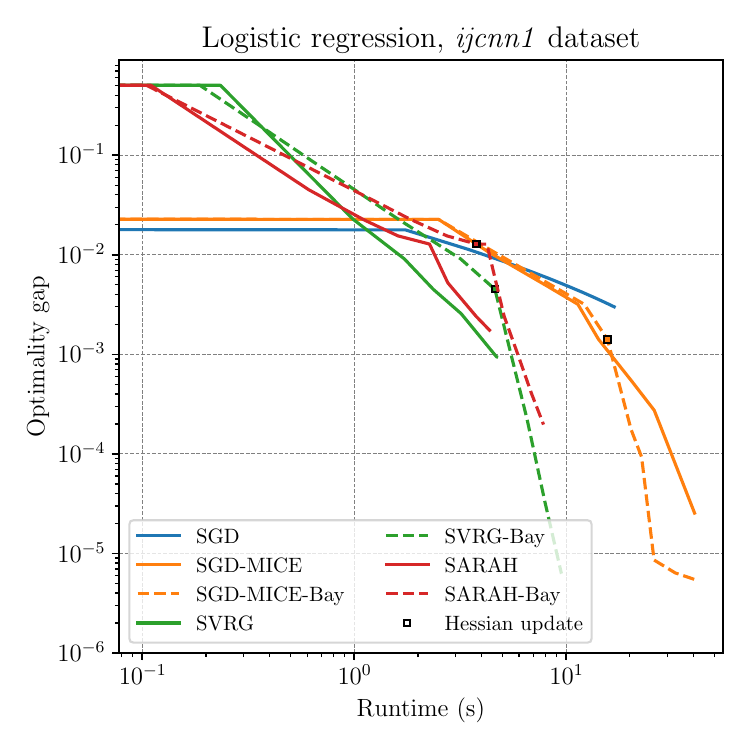}
 \caption{Convergence of the optimality gap versus iteration (left) and runtime (right) for Example \ref{sec:logreg}, \emph{ijcnn1} dataset. The Hessian updates are marked with colored squares with black edges.}
 \label{fig:ex2_ijcnn1_conv}
\end{figure}

\begin{figure}
  \centering
 \includegraphics[width=.6\linewidth]{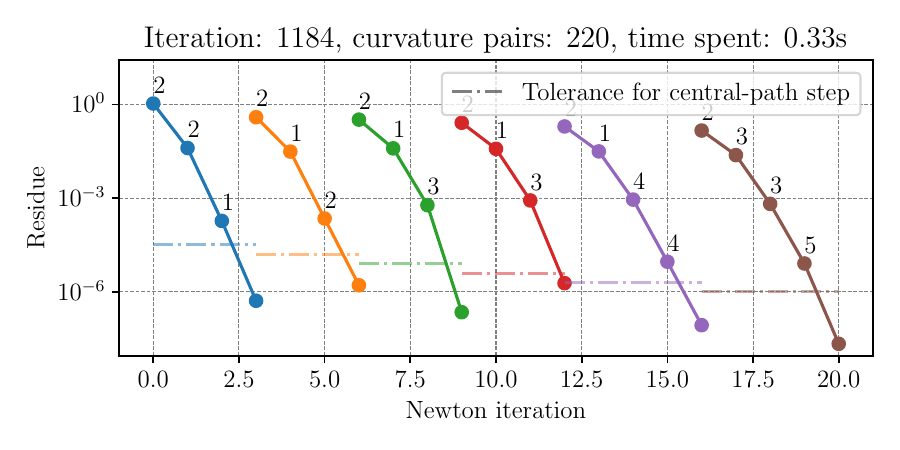}
 \caption{Convergence of the Newton-CG method for each of the Hessian updates of SGD-MICE-Bay in Example \ref{sec:logreg}, \emph{ijcnn1} dataset.
 The number of curvature pairs available and the time taken are presented for each Hessian update.
 Each color represents a different step of the central-path method with a different logarithmic barrier parameter $\beta$ and a different residue tolerance.
 The tolerance for each central-path step is presented as a dash-dotted line.
 The number above each Newton iteration represents the number of conjugate gradient iterations needed to find the Newton direction.
 }
 \label{fig:ex2_ijcnn1_newton_cg}
\end{figure}

Figure \ref{fig:ex2_ijcnn1_conv} presents the convergence for the \emph{ijcnn1} dataset.
Here, as in the case of the \emph{mushrooms} dataset, our Bayesian approach improved the convergence in all cases.
However, the most remarkable improvement is the one of SVRG-Bay, which achieved an optimality gap of more than two orders of magnitude smaller than the vanilla SVRG with just one Hessian update.
Both SVRG-Bay and SGD-MICE-Bay achieved very similar values of optimality gap, with SVRG-Bay requiring less runtime than SGD-MICE-Bay.
In Figure \ref{fig:ex2_ijcnn1_newton_cg}, we present the convergence of Newton-CG for the one Hessian update of SGD-MICE-Bay for the \emph{ijcnn1} dataset.
Note that the time spent computing the Hessian approximation is 0.33 seconds; for comparison, computing the true Hessian took 0.54 seconds in an average of 100 evaluations.

\begin{figure}
  \includegraphics[width=.48\linewidth]{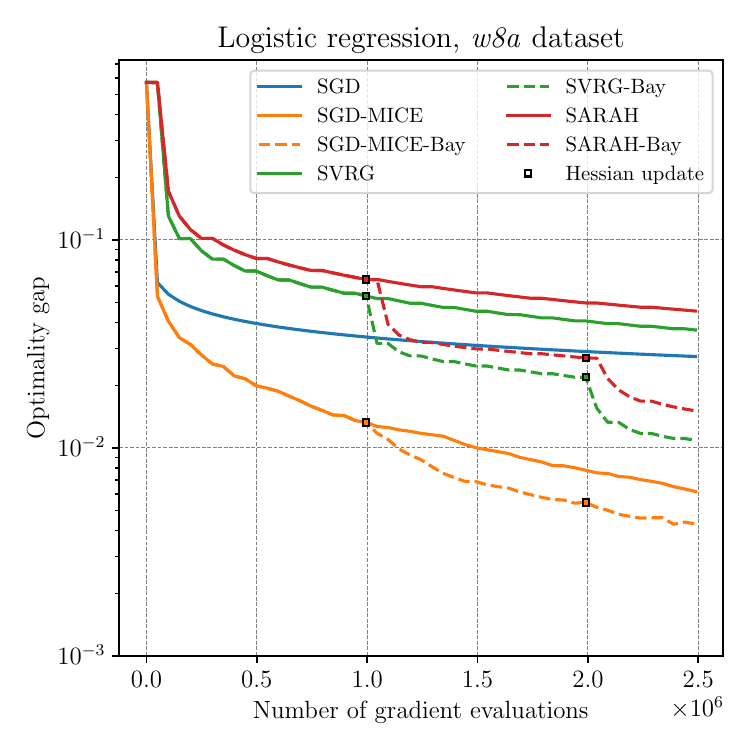}
  \includegraphics[width=.48\linewidth]{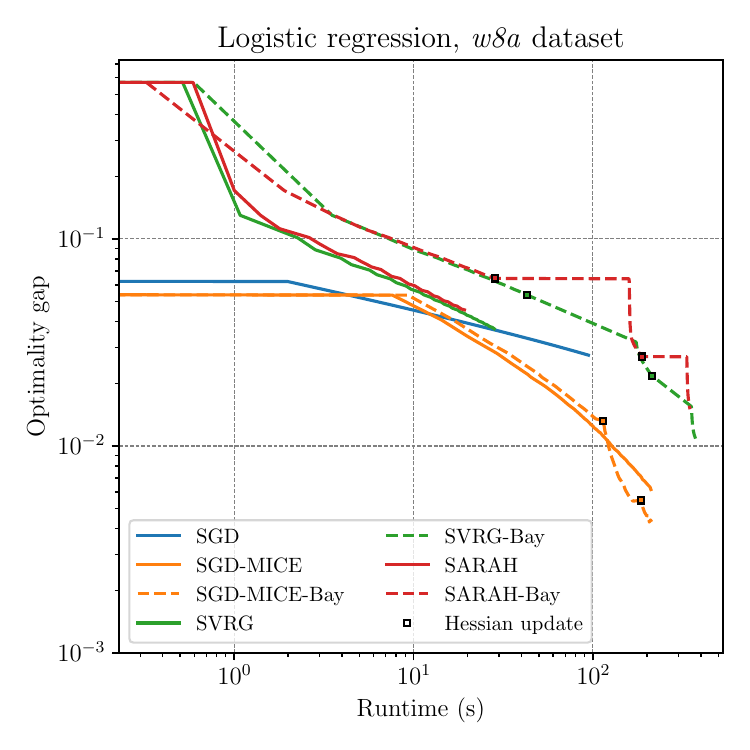}
 \caption{Convergence of the optimality gap versus iteration (left) and runtime (right) for Example \ref{sec:logreg}, \emph{w8a} dataset. The Hessian updates are marked with colored squares with black edges.}
 \label{fig:ex2_w8a_conv}
\end{figure}

\begin{figure}
  \centering
 \includegraphics[width=.6\linewidth]{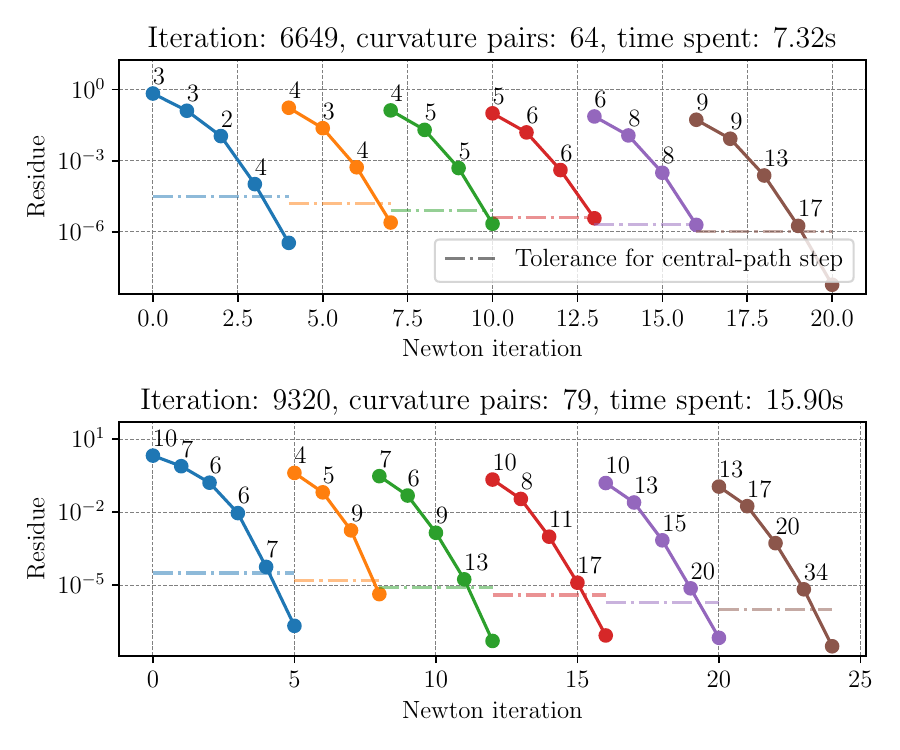}
 \caption{Convergence of the Newton-CG method for each of the Hessian updates of SGD-MICE-Bay in Example \ref{sec:logreg}, \emph{w8a} dataset.
 The number of curvature pairs available and the time taken are presented for each Hessian update.
 Each color represents a different step of the central-path method with a different logarithmic barrier parameter $\beta$ and a different residue tolerance.
 The tolerance for each central-path step is presented as a dash-dotted line.
 The number above each Newton iteration represents the number of conjugate gradient iterations needed to find the Newton direction.
 }
 \label{fig:ex2_w8a_newton_cg}
\end{figure}

The \emph{w8a} dataset is the case with the largest number of features, 300, meaning that the Hessian to be approximated by our Bayesian approach has $300 \times 300$ components.
The convergence for the \emph{w8a} dataset is presented in Figure \ref{fig:ex2_w8a_conv}.
Coupling SVRG, SARAH, and SGD-MICE with our Bayesian approach improved the convergence in all cases, but it was more significant for SVRG and SARAH.
However, in terms of runtime, SGD-MICE was able to perform better than its competitors, and SGD-MICE-Bay had a small edge over SGD-MICE.
Figure \ref{fig:ex2_w8a_newton_cg} presents the convergence of Newton-CG for SGD-MICE-Bay on the \emph{w8a} dataset.
As in the case of the \emph{mushrooms} dataset, a small number of Newton steps (a maximum of $24$) and conjugate gradient steps (a maximum of $34$) were needed despite the large dimensionality of the problem at $90000$ variables.
Even with the large Hessian, the time spent to find an approximation ($7.32$ s in the first update and $15.9$ s in the second) is competitive with the time of computing the true Hessian, of $8.52$ s in an average of 100 evaluations.
Remember that our approach has stability advantages over using the true Hessian inverse as a preconditioner due to the control of the noise amplification given by the logarithmic barrier on the smallest eigenvalue.

\begin{figure}
  \centering
  \includegraphics[width=.48\linewidth]{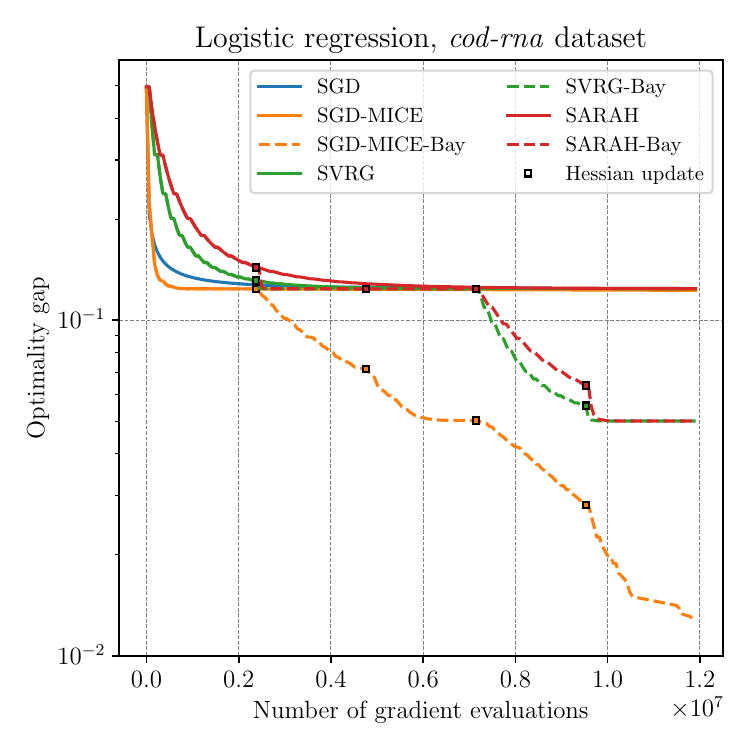}
  \includegraphics[width=.48\linewidth]{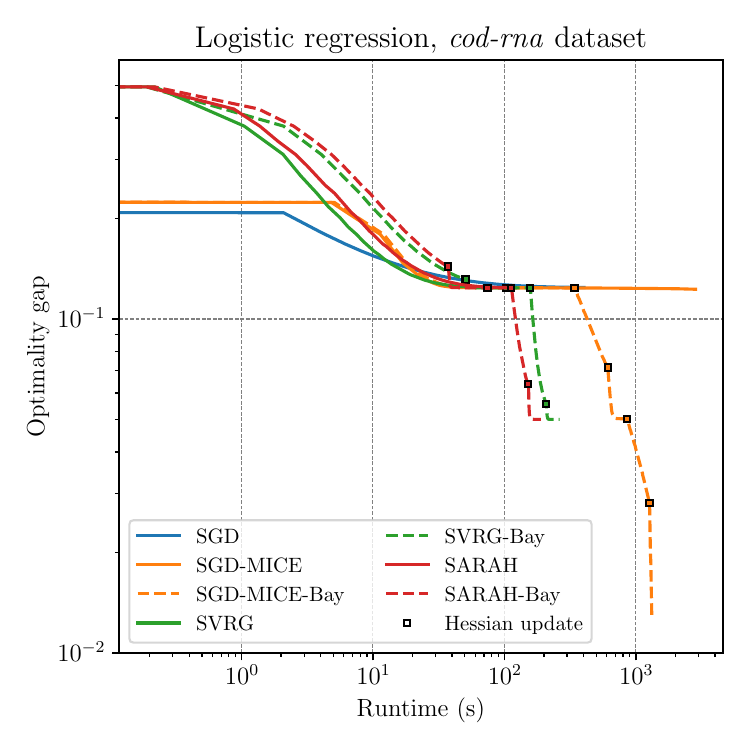}
 \caption{Convergence of the optimality gap versus iteration (left) and runtime (right) for Example \ref{sec:logreg}, \emph{cod-rna} dataset. The Hessian updates are marked with colored squares with black edges.}
 \label{fig:ex2_cod_rna_conv}
\end{figure}

\begin{figure}
  \centering
 \includegraphics[width=.6\linewidth]{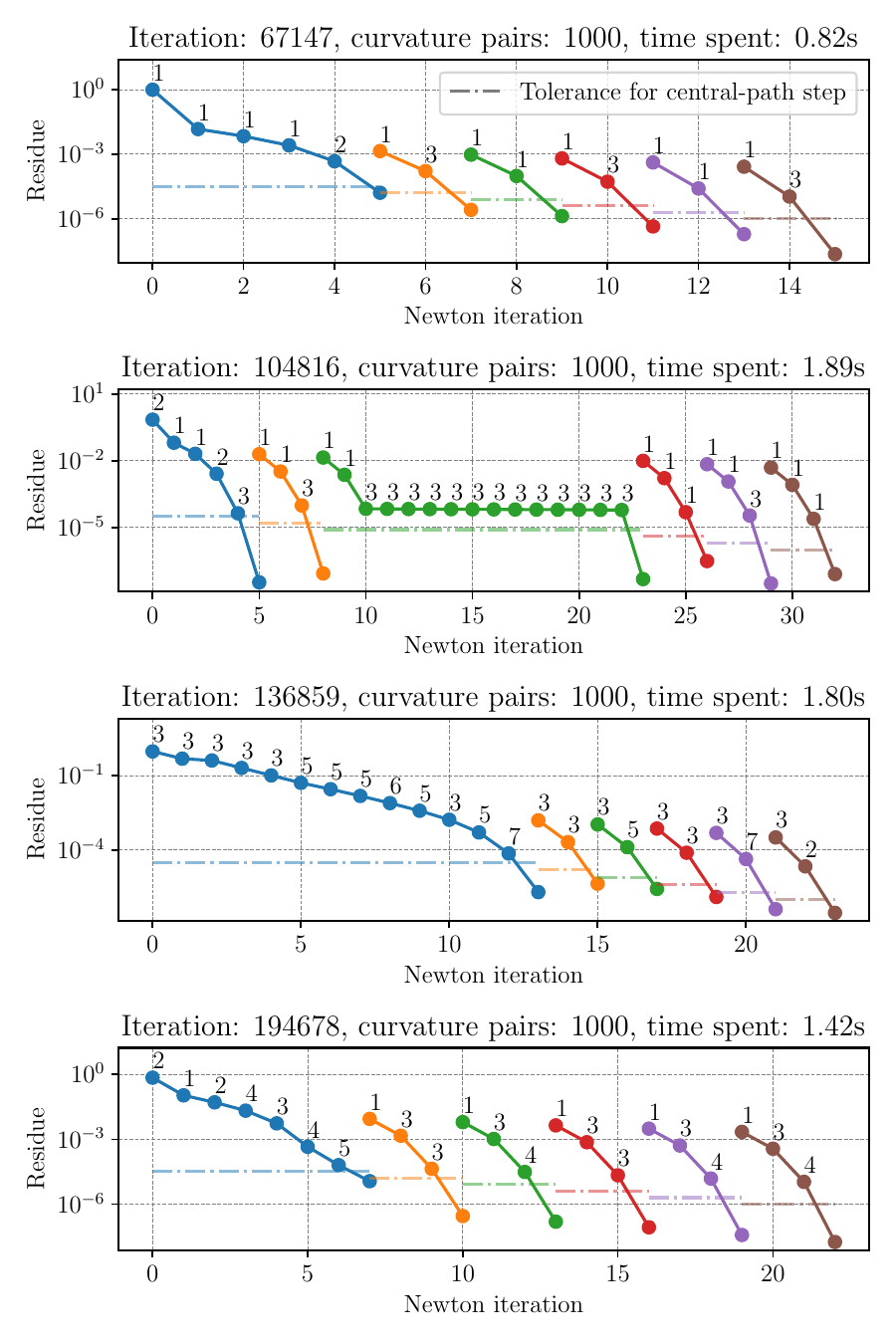}
 \caption{Convergence of the Newton-CG method for each of the Hessian updates of SGD-MICE-Bay in Example \ref{sec:logreg}, \emph{cod-rna} dataset.
 The number of curvature pairs available and the time taken are presented for each Hessian update.
 Each color represents a different step of the central-path method with a different logarithmic barrier parameter $\beta$ and a different residue tolerance.
 The tolerance for each central-path step is presented as a dash-dotted line.
 The number above each Newton iteration represents the number of conjugate gradient iterations needed to find the Newton direction.
 }
 \label{fig:ex2_cod_rna_newton_cg}
\end{figure}

Training the logistic regression model to the \emph{cod-rna} dataset is the problem with the largest condition number of those studied here, at $6.48 \times 10^9$.
We expect our Bayesian approach to be able to assist the optimization methods in converging to the true solution despite the large condition number of the problem.
The convergence of the optimality gap for the \emph{cod-rna} dataset is presented in Figure \ref{fig:ex2_cod_rna_conv}.
Note that SVRG, SARAH, and SGD-MICE get stuck after some iterations, improving very little for most of the optimization process.
Their preconditioned counterparts, however, perform much better, being able to decrease the optimality gap continuously.
When looking at the convergence of the optimality gap versus runtime for SGD-MICE-Bay, we observe a sharp decrease after the first Hessian update that follows until its stop.
Figure \ref{fig:ex2_cod_rna_newton_cg} presents the convergence of Newton-CG for the four Hessian updates of SGD-MICE-Bay.
In the second Hessian update, the third central-path step, the Newton-CG method got temporarily stuck for 12 iterations, leading us to believe that the parameters $\beta$ and/or $\rho$ were not properly set for this example.
Still, the Newton-CG method could resume the Hessian approximation in a reasonable time, at $1.89$ s, showing the robustness of our approach.

\begin{figure}
  \includegraphics[width=.48\linewidth]{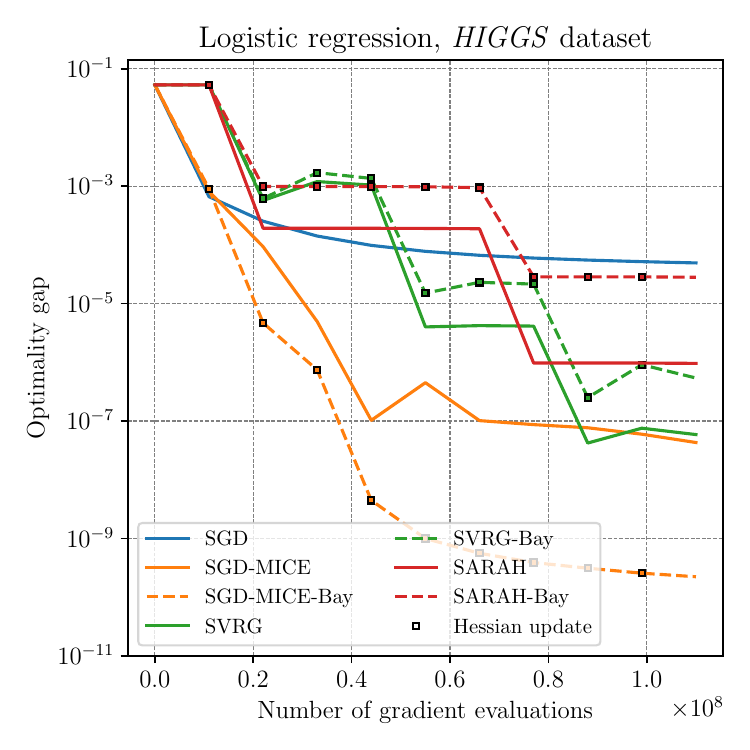}
  \includegraphics[width=.48\linewidth]{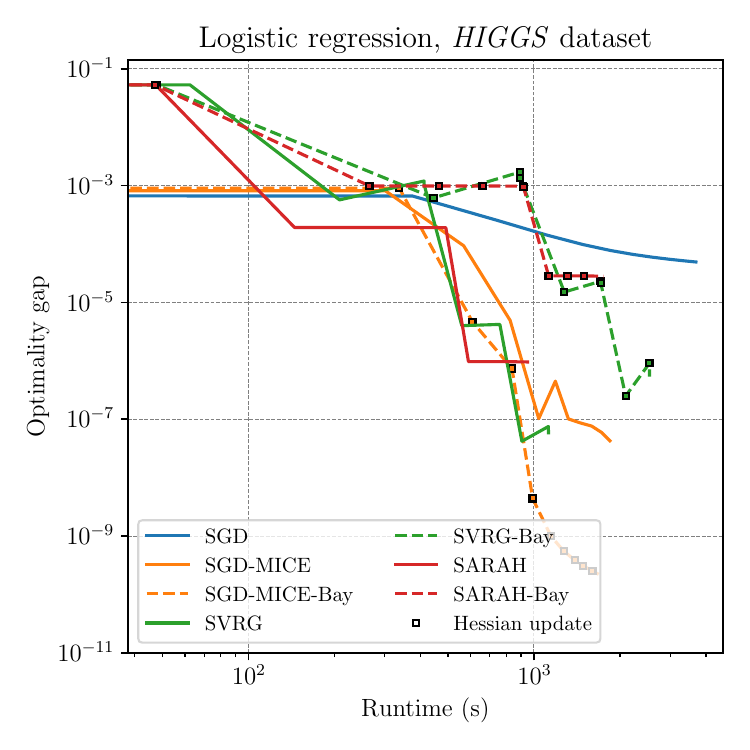}
 \caption{Convergence of the optimality gap versus iteration (left) and runtime (right) for Example \ref{sec:logreg}, \emph{HIGGS} dataset. The Hessian updates are marked with colored squares with black edges.}
 \label{fig:ex2_higgs_conv}
\end{figure}

The \emph{HIGGS} dataset contains $11 \times 10^6$ data points, being the largest studied here.
Thus, since our Hessian approximation is independent of the data size, we expect to have a time advantage compared to computing the true Hessian.
Figure \ref{fig:ex2_higgs_conv} presents the convergence of the optimality gap for the \emph{HIGGS} dataset.
Here we observe that, in contrast to the previous results, our Bayesian approach worsened the performance of SARAH and SVRG.
For SGD-MICE, however, our method improved the final optimality gap by more than one order of magnitude.
Also, concerning runtime, SGD-MICE-Bay performed better than the other methods, despite doing 9 Hessian updates.
The longest Hessian update was the last one, lasting $2.36$ s.
In comparison, the average time taken to compute a Hessian for this example is $130.58$ s.
One possible interpretation for the bad performance of SVRG-Bay and SARAH-Bay is that, contrary to SGD-MICE, SVRG and SARAH do not control the relative statistical error in the gradient estimates.

\section{Acknowledgments}

This work was partially supported by the KAUST Office of Sponsored Research (OSR) under Award numbers URF$/1/2281-01-01$, URF$/1/2584-01-01$ in the KAUST Competitive Research Grants Program Round 8, the Alexander von Humboldt Foundation.

\section{Conclusion}

We presented a Bayesian approach to approximate Hessians of functions given noisy observations of their gradients.
This problem is of great importance in stochastic optimization, where pre-conditioning gradient estimates with the inverse of the Hessian matrix can improve the convergence of stochastic gradient descent methods.
The proposed Bayesian approach takes into consideration not only the secant equations as in deterministic quasi-Newton methods but also the noise in the gradients.
To mitigate the known effect of amplification of the gradient noise, we control the smallest eigenvalue of the Hessian approximation.
To maximize the log-posterior of the Hessian matrix, we use a Newton-CG method with a central-path approach to impose the eigenvalue constraints.
The numerical results presented show that our approach is not only interesting from the theoretical perspective but also results in practical advantage.
In both a stochastic quadratic equation and in the training of a logistic regression model with $\ell_2$ regularization, using our Bayesian approximation of the Hessian matrix resulted in better convergence rates without a significant increase in runtime.
For future research, we suggest using low-rank or diagonal approximations of the Hessians.




\end{document}